\documentclass[10pt]{article}
\usepackage{bbm}
\usepackage{amssymb}
\usepackage{mathrsfs}
\usepackage{calc}
\usepackage{graphicx}
\usepackage[all,cmtip]{xy}
\usepackage{amsmath}
\usepackage{amscd}
\usepackage{lmodern}
\usepackage{amsfonts}
\usepackage{picinpar}
\usepackage{tikz-cd}
\usepackage{pgf,tikz,pgfplots}
\pgfplotsset{compat=1.14}
\usepackage{mathrsfs}
\usetikzlibrary{arrows}
\usepackage[amsmath,thmmarks]{ntheorem}
\usepackage{enumerate}
\usepackage{mathptmx}
\usepackage{stmaryrd}
\usepackage{array}
\usepackage{extarrows}
\usepackage[nottoc,numbib]{tocbibind}
\usepackage[]{hyperref}

\newtheorem{defi}{Definition}[section]
\newtheorem{thm}[defi]{Theorem}

\newtheorem{prop}[defi]{Proposition}
\newtheorem{cor}[defi]{Corollary}
\newtheorem{lem}[defi]{Lemma}
\theorembodyfont{}
\newtheorem{exa}[defi]{Example}
\theoremstyle{plain}
\newtheorem{rmk}[defi]{Remark}
\newtheorem{question}[defi]{Question}

\theoremstyle{nonumberplain} \theorembodyfont{} \theoremsymbol{\ensuremath{\square}} \qedsymbol{\ensuremath\square}
\newtheorem{pro}{Proof}

\DeclareMathOperator{\GL}{GL}

\DeclareMathOperator{\PGL}{PGL}

\DeclareMathOperator{\Aff}{Aff}
\DeclareMathOperator{\Bir}{Bir}
\DeclareMathOperator{\Aut}{Aut}

\DeclareMathOperator{\Id}{Id}

\DeclareMathOperator{\Ind}{Ind}

\begin{document}

\title{Uniqueness of birational structures on Inoue surfaces}
\author{ZHAO ShengYuan}\date{}
\maketitle

\newcommand{\CRC}{\operatorname{Bir}(\mathbb{P}^2)}
\newcommand{\CRn}{\operatorname{Cr}_{n}(\mathbf{C})}
\newcommand{\Jonq}{\operatorname{Jonq}(\mathbf{C})}
\newcommand{\Jonqz}{\operatorname{Jonq}_0(\mathbf{C})}
\newcommand{\Dev}{\operatorname{Dev}}
\newcommand{\Hol}{\operatorname{Hol}}
\newcommand{\CC}{\mathbf{C}}
\newcommand{\RR}{\mathbf{R}}
\newcommand{\NN}{\mathbf{N}}
\newcommand{\QQ}{\mathbf{Q}}
\newcommand{\ZZ}{\mathbf{Z}}
\newcommand{\FF}{\mathscr{F}}
\newcommand{\FFh}{\mathbb{F}}
\newcommand{\KK}{\mathbf{K}}
\newcommand{\Sym}{\mathscr{S}}
\newcommand{\HH}{\mathbb{H}}
\newcommand{\PP}{\mathbb{P}}
\newcommand{\DD}{\mathbb{D}}
\newcommand{\Jonqui}{Jonqui\`eres }
\newcommand{\Kah}{K\"{a}hler}

\newenvironment{ital}{\it}{}

\begin{abstract}
We prove that the natural $(\Aff _2(\CC),\CC^2)$-structure on an Inoue surface is the unique $(\CRC,\PP^2(\CC))$-structure, generalizing a result of Bruno Klingler which asserts that the natural $(\Aff _2(\CC),\CC^2)$-structure is the unique $(\PGL _3(\CC),\PP^2(\CC))$-structure.
\end{abstract}
\renewcommand{\abstractname}{Abstract}

\setcounter{tocdepth}{2}
\tableofcontents

\section{Introduction}
Inoue surfaces are compact non-\Kah{} complex surfaces discovered by Inoue in \cite{Ino74}. They are of class VII in Enriques-Kodaira's classification of compact complex surfaces, and are the only compact complex sufaces with Betti numbers $b_1=1,b_2=0$ (cf. \cite{Tel94}). There are three types of Inoue surfaces: $S^0$, $S^+$ and $S^{-}$. Their universal covers are isomorphic to $\HH\times\CC$ where $\HH$ is the upper half plane and the deck transformations can be written as restrictions of complex affine transformations of $\CC^2$. Therefore the Inoue surfaces are equipped with a natural complex affine structure. Klingler proves in \cite{Kli98} that the natural complex affine structure is the unique complex projective structure carried by Inoue surfaces. In this article we prove the following:
\begin{thm}\label{mainthm}
If an Inoue surface is equipped with a $(\Bir(X),X)$-structure for some complex projective surface $X$, then $X$ is a rational surface and the $(\Bir(X),X)$-structure is induced by the natural $(\Aff _2(\CC),\CC^2)$-structure.
\end{thm}

\begin{rmk}
Roughly speaking, a birational structure is an atlas of local charts with birational changes of coordinates. The precise definition and basic properties will be given in Section \ref{birationalstructure}. It is a generalization of the classical $(G,X)$-structure; if we think of a geometric structure as a way to patch coordinates, then it is the most general algebraic geometric structure (the changes of coordinates are rational). 

In a recent article \cite{KS19}, Kwon and Sullivan introduced some generalized notions of geometric structures for which they allow a family of Lie groups $(G_i)_{i}$ acting on the same space $X$. The group of birational transformations of a surface, though not a classical Lie group itself, is generated by Lie groups acting by holomorphic diffeomorphisms on different birational models of $X$. So, the geometric structure of \cite{KS19} share interesting similarities with $(\Bir(X),X)$-structures. Note that the group of birational transformations of a variety of dimension $\geq 3$ may not be generated by its connected algebraic subgroups (cf. \cite{BY19}). Kwon and Sullivan proved in \cite{KS19} that every prime orientable three manifold admits such a generalized geometric structure. Their result is somewhat analogous to Dloussky's conjecture mentionned in Remark \ref{Dloussky} below.
\end{rmk}

\begin{rmk}
Compared to the four-pages-long proof in \cite{Kli98} of the uniqueness of complex projective structure, our proof is more involved because the group $\CRC$ of birational transformations of $\PP^2$ is much larger than $\PGL_3(\CC)$. Also the fact that in our case the developing map is a priori not holomorphic but only meromorphic will be the cause of some technical complications.
\end{rmk}

\begin{rmk}\label{Dloussky}
Complex projective structures on compact complex surfaces are classified by Klingler in \cite{Kli98}. There exist compact non-\Kah{} complex surfaces which have $(\CRC,\PP^2)$-structures but no complex projective structures, for example some Hopf surfaces; and Dloussky conjectured in \cite{Dlo16} that every surface of class VII admits a $(\CRC,\PP^2)$-structure. 
\end{rmk}

\begin{rmk}
If $Y$ is a complex projective surface, we say that a subgroup $\Gamma$ of $\Bir(Y)$ has the \emph{Kleinian property} if the following three conditions are satisfied: 1) the group $\Gamma\subset\Bir(Y)$ acts by holomorphic diffeomorphisms on a Euclidean open set $U\subset Y$, i.e.\ an open set for the Euclidean topology but not necessarily for the Zariski topology; 2) the action of $\Gamma$ on $U$ is free and properly discontinuous; 3) the quotient $U/\Gamma$ is compact. Once we have a birational Kleinian group, the quotient surface is equipped naturally with a birational structure. Thus, we can view Theorem \ref{mainthm} as a result about Fatou components of (groups of) birational transformations. For a systematic study of birational Kleinian groups we refer to the forthcoming article \cite{Zhaoklein}. 
\end{rmk}

\paragraph{Plan and strategy.}
Section \ref{firstsection} concerns two subjects of independent interest. The notion of birational structure appeared already in the work of Dloussky \cite{Dlo16} but several subtleties, that do not appear in classical geometric structures, were not addressed in that paper. In Section \ref{birationalstructure} we give two different definitions of $(\Bir(\PP^n),\PP^n)$-structure. For $n=2$ they are the same but for $n\geq 3$ whether they are the same is equivalent to an open question of Gromov. In Section \ref{AhlforsNevanlinnasection} we study Ahlfors-Nevanlinna currents attached to entire curves (see the work of Brunella and McQuillan in \cite{Mcq98}, \cite{Bru99}): for nice (uniform) families of entire curves, we show how to construct families of Ahlfors-Nevalinna 
currents, with a fixed cohomology class; this may be useful to people interested in holomorphic foliations or Kobayashi
hyperbolicity. 

After these preliminaries we prove Theorem \ref{mainthm}. The construction of Inoue surfaces of type $S^0$ (resp. $S^{\pm}$) will be recalled in Section \ref{secondsection} (resp. \ref{thirdsection}). For simplicity, let us focus, here, on Inoue surfaces of type $S^0$. Rather different tools are used, depending on the size of the image of the holonomy representation. When the holonomy is injective, the two principal ingredients are the classification of solvable subgroups of $\CRC$ due to Cantat (\cite{Can11}), D\'eserti (\cite{Des15}) and Urech (\cite{Ure}), and the classification of subgroups of $\CRC$ ismorphic to $\ZZ^2$, obtained in \cite{Zhaocent}. With these results, we can reduce the structure group from $\CRC$ to $\PGL_3(\CC)$, and then apply Klingler's previous work \cite{Kli98}. When the holonomy representation is not injective, we can suppose that its image is cyclic. Then, the strategy is geometric: an Inoue surface is foliated by compact real submanifolds of dimension three that are themselves foliated by entire curves, i.e.\ Riemann surfaces isomorphic to $\CC$. Via the developing map, we obtain families of Levi-flat hypersurfaces foliated by entire curves, in some projective surfaces. The proof is then based on the following three tools that are described in Section \ref{firstsection}: 1) our deformation lemma for Ahlfors-Nevanlinna currents; 2) the relation between these currents and the transverse invariant measures of Plante (\cite{Pla75}) and Sullivan (\cite{Sul76}); 3) properties of the pull-back action of a birational transformation on currents (as in \cite{DF01} and \cite{Can01}). 

\paragraph{Acknowledgement.}
It is a pleasure to thank Serge Cantat for his constant support and numerous discussions. I would also like to thank J\'er\'emy Blanc, Bertrand Deroin and JunYi Xie for interesting discussions, B. D. for the reference \cite{KS19} and J. B. for pointing out to me Question \ref{gromovquestion}. My thanks also goes to the referee whose comments have been helpful to improve the text.

\section{Preliminaries}\label{firstsection}
\subsection{Groups of birational transformations}\label{birpre}
Let $X$ be a smooth complex projective surface. We denote by $\Bir(X)$ the group of birational transformations of $X$. An element $f$ of $\Bir(X)$ has a pull back action $f^*$ on $\operatorname{H}^{1,1}(X,\RR)$ (cf. \cite{DF01}). Note that in general $(f^*)^n\neq(f^n)^*$. Fix an ample class $H\in \operatorname{H}^{1,1}(X,\RR)$, the \emph{$H$-degree} of $f$ is the intersection number $f^*H\cdot H$. 

The \emph{plane Cremona group} $\CRC$ is the group of birational transformations of the projective plane $\PP^2_{\CC}$. It is isomorphic to the group of $\CC$-algebra automorphisms of $\CC(X_1,X_2)$, the function field of $\PP^2_{\CC}$. Using a system of homogeneous coordinates $[x_0;x_1;x_2]$, a birational transformation $f\in\CRC$ can be written as
\[
[x_0:x_1:x_2]\dashrightarrow [f_0(x_0,x_1,x_2):f_1(x_0,x_1,x_2):f_2(x_0,x_1,x_2)]
\]
where $f_0,f_1,f_2$ are homogeneous polynomials of the same degree without common factor. This degree does not depend on the system of homogeneous coordinates and is the degree of $f$ with respect to the class of a projective line. Birational transformations of degree $1$ are homographies and form $\Aut (\PP^2) = \PGL _3 (\CC)$, the group of automorphisms of the projective plane. See \cite{Can18} for more about the Cremona group.

\paragraph{Algebraically stable maps.}
If $f$ is a birational transformation of a smooth projective surface $X$, we denote by $\Ind(f)$ the set of indeterminacy points of $f$. We say that $f$ is \emph{algebraically stable} if there are no curves $V$ on $X$ such that $f^k(V)\subset \Ind(f)$ for some integer $k\geq 0$. There always exists a birational morphism $\hat{X}\rightarrow X$ which lifts $f$ to an algebraically stable birational transformation of $\hat{X}$ (\cite{DF01} Theorem 0.1). An algebraically stable map $f$ satisfies $(f^*)^n=(f^n)^*$ (cf. \cite{DF01}). 

\paragraph{Four types of elements.}
 Fix a Euclidean norm $\left\|\cdot\right\|$ on $\operatorname{H}^{1,1}(X,\RR)$. The two sequences $\left(\left\|(f^n)^*\right\|\right)_n$ and $((f^n)^*H\cdot H)_n$ have the same asymptotic growth. Elements of $\Bir(X)$ are classified into four types (cf. \cite{DF01}):  
\begin{enumerate}
	\item The sequence $\left(\left\|(f^n)^*\right\|\right)_{n\in \NN}$ is bounded, $f$ is birationally conjugate to an automorphism of a smooth birational model of $X$ and a positive iterate of $f$ lies in the connected component of identity of the automorphism group of that surface. We call $f$ an \emph{elliptic} element.
	\item The sequence $\left(\left\|(f^n)^*\right\|\right)_{n\in \NN}$ grows linearly, $f$ preserves a unique pencil of rational curves and $f$ is not conjugate to an automorphism of any birational model of $X$. We call $f$ a \emph{\Jonqui twist}.
	\item The sequence $\left(\left\|(f^n)^*\right\|\right)_{n\in \NN}$ grows quadratically, $f$ is conjugate to an automorphism of a smooth birational model preserving a unique genus one fibration. We call $f$ a \emph{Halphen twist}.
	\item The sequence $\left(\left\|(f^n)^*\right\|\right)_{n\in \NN}$ grows exponentially and $f$ is called \emph{loxodromic}. The limit $\lambda(f)=\lim_{n\rightarrow +\infty}\left(\left\|(f^n)^*\right\|\right)^{\frac{1}{n}}$ exists and we call it the \emph{dynamical degree} of $f$. If $f$ is an algebraically stable map on $X$, then there is a nef cohomology class $v_f^+\in \operatorname{H}^{1,1}(X,\RR)$, unique up to multiplication by a constant, such that $f^*v_f^+=\lambda(f)v_f^+$. If moreover $v_f^+$ has zero self-intersection, then $f$ is conjugate to an automorphism.
\end{enumerate}

\paragraph{Loxodromic automorphisms.}
We refer the reader to \cite{Can14} for details of the materials presented in this paragraph. Let $X$ be a smooth projective surface and $f$ be an automorphism of $X$ which is loxodromic. The dynamical degree $\lambda(f)$ is a simple eigenvalue for the pullback action $f^*$ on $\operatorname{H}^{1,1}(X,\RR)$ and it is the unique eigenvalue of modulus larger than $1$. Let $v_f^+\in \operatorname{H}^{1,1}(X,\RR)$ be a non-zero eigenvector associated with $\lambda(f)$; we have $f^*v_f^+=\lambda(f)v_f^+$. By considering $f^{-1}$, we can also find a non-zero eigenvector $v_f^-$ such that $f^*v_f^-=\frac{1}{\lambda(f)}v_f^-$. The two cohomology classes $v_f^+,v_f^-$ are nef and of self-intersection $0$, they are uniquely determined up to scalar multiplication. They are irrational in the sense that the two lines $\RR v_f^+$ and $\RR v_f^-$ contain no non-zero elements of $\operatorname{H}^{1,1}(X,\RR)\cap \operatorname{H}^{2}(X,\ZZ)$. We will need the following theorem of Cantat which has been generalized to higher dimension by Dinh and Sibony:

\begin{thm}[Cantat \cite{Can01}, \cite{Moncet}, \cite{DS05}, \cite{DS10}]\label{CantatDinhSibony}
There is a unique closed positive current $T_f^+$ (resp. $T_f^-$) whose cohomology class is $v_f^+$ (resp. $v_f^-$). It satisfies $f^*T_f^+=\lambda(f)T_f^+$ (resp. $f^*T_f^-=\frac{1}{\lambda(f)}T_f^-$). 
\end{thm}
   
\paragraph{The \Jonqui group.}
Fix an affine chart of $\PP^2$ with coordinates $(x,y)$.
\emph{The \Jonqui group} $\Jonq$ is the subgroup of the Cremona group of all transformations of the form
\[
(x,y)\mapsto \left ( \frac{ax+b}{cx+d},\frac{A(x)y+B(x)}{C(x)y+D(x)} \right ),\quad 
\begin{pmatrix}
	a&b\\c&d
\end{pmatrix}\in \PGL _2 (\CC),\quad 
\begin{pmatrix}
	A&B\\C&D
\end{pmatrix}\in \PGL _2(\CC(x)).
\]
In other words, the \Jonqui group is the maximal group of birational transformations of $\PP^1\times\PP^1$ permuting the fibers of the
projection onto the first factor; it is isomorphic to the semidirect product $\PGL _2(\CC)\ltimes \operatorname{PGL}_2(\CC(x))$. A different choice of the affine chart yields a conjugation by an element of $\PGL _3(\CC)$. More generally a conjugation by an element of the Cremona group yields a maximal group preserving a pencil of rational curves; conversely any two such groups are conjugate in $\CRC$. 

Elements of the \Jonqui group are either elliptic or \Jonqui twists. We will need the following results:
\begin{thm}[\cite{Zhaocent}]\label{centralizer}
Let $G$ be a subgroup of $\Jonq$ which is isomorphic to $\ZZ^2$. Then $G$ has a pair of generators $(f,g)$ such that one of the following (mutually exclusive) situations happens:
\begin{enumerate}
	\item $f,g$ are elliptic elements and $G\subset \Aut(X)$ where $X$ is a rational surface; 
	\item $f$ is a \Jonqui twist, and a finite index subgroup of $G$ preserves each fiber of the $f$-invariant fibration; 
	\item $f$ is a \Jonqui twist with action of infinite order on the base of the rational fibration and $g$ is an elliptic element whose action on the base is of finite order. In some affine chart, we can write $f,g$ in one of the following forms:
	\begin{itemize}
		\item $g$ is $(x,y)\mapsto(\alpha x,\beta y)$ and $f$ is $(x,y)\dashrightarrow(\eta(x),yR(x^k))$ where $\alpha,\beta\in\CC^*, \alpha^k=1, R\in \CC(x),\eta\in \PGL_2(\CC), \eta(\alpha x)=\alpha \eta(x)$ and $\eta$ is of infinite order;
		\item $g$ is $(x,y)\mapsto(\alpha x,y+1)$ and $f$ is $(x,y)\dashrightarrow(\eta(x),y+R(x))$ where $\alpha\in\CC^*, R\in \CC(x), R(\alpha x)=R(x), \eta\in \PGL_2(\CC), \eta(\alpha x)=\alpha \eta(x)$ and $\eta$ is of infinite order.
	\end{itemize}
\end{enumerate}
\end{thm}

\begin{thm}[\cite{BD15} Lemmata 2.7 and 2.8]\label{ellipticthm}
Let $f\in\CRC$ be an elliptic element of infinite order. 
\begin{enumerate}
	\item If $f$ is of the form $(x,y)\mapsto(x,\nu y)$ where $\nu \in \CC^*$ has infinite order, then the centralizer of $f$ in $\CRC$ is
	\[
	\{(x,y)\dashrightarrow(\eta(x),yR(x))\vert \eta\in \PGL_2(\CC),R\in \CC(x)\}.
	\]
	\item If $f$ is of the form $(x,y)\mapsto(x,y+v)$ with $v\in\CC^*$, then the centralizer of $f$ in $\CRC$ is
	\[
  \{(x,y)\dashrightarrow(\eta(x),y+R(x))\vert \eta\in \PGL_2(\CC), R\in \CC(x)\}.
	\]
\end{enumerate}
\end{thm}

\paragraph{Tits alternative and solvable subgroups.}
D\'eserti and Urech refined for finitely generated solvable subgroups, the strong Tits alternative proved by Cantat in \cite{Can11}; we state the solvable version:
\begin{thm}[Cantat, D\'eserti, Urech \cite{Can11}, \cite{Des15}, \cite{Ure}]\label{Deserti}
Let $G\subset\Bir(X)$ be a solvable subgroup. Exactly one of the following cases holds up to conjugation.
  \begin{enumerate}
  \item $G$ is a subgroup of automorphisms of a birational model $Y$ and a finite index subgroup of $G$ is in $\Aut(Y)^0$ the connected component of identity of $\Aut(Y)$; the elements of $G$ are all elliptic and $G$ is called an elliptic subgroup.
  \item $G$ preserves a rational fibration and has at least one \Jonqui twist. 
  \item $G$ is a virtually abelian group whose elements are Halphen twists; there is a birational model $Y$ on which the action of $G$ is by automorphisms and preserves an elliptic fibration. 
  \item $X$ is a rational surface and $G$ is contained in the group generated by $\{(\alpha x,\beta y)\vert \alpha,\beta \in\CC^*\}$ and one loxodromic monomial transformation $(x^py^q,x^ry^s)$ where $\begin{pmatrix}p&q\\r&s\end{pmatrix}\in \operatorname{GL}_2(\ZZ)$ is a hyperbolic matrix.
	\item $X$ is an abelian surface and $G$ is contained in the group generated by translations and one loxodromic transformation.
 \end{enumerate}
\end{thm}

\subsection{Geometric structures}\label{birationalstructure}
Let us first recall the classical notion of $(G,X)$-structures in the sense of Ehresmann (cf. \cite{Ehr36}, see also \cite{Thu97}):
\begin{defi}
Let $X$ be a connected real analytic manifold and let $G$ be a Lie group which acts real analytically faithfully on $X$. Let $V$ be a real analytic manifold. A $(G,X)$-structure on $V$ is a maximal atlas of local charts $\phi_i:U_i\rightarrow X$ such that
\begin{itemize}
	\item the $U_i$ are open sets of $V$ and form a covering;
	\item the $\phi_i$ are diffeomorphisms onto their images;
	\item the changes of coordinates $\phi_i \circ\phi_j^{-1}:\phi_j(U_i\cap U_j)\rightarrow\phi_i(U_i\cap U_j)$ are restrictions of elements of $G$.
\end{itemize}
A $(G,X)$-manifold is a manifold which is equipped with a $(G,X)$-structure.
\end{defi}

The group of birational transformations of an algebraic variety is not a Lie group in the classical sense, see \cite{BF13} for its topology. Its action on the variety is not a classical set-theoretic group action either. We give here two non-equivalent definitions of birational structure. The first one is more flexible and is the notion of birational structure that we use in this article. 
\begin{defi}
Let $V$ be a complex manifold. Let $X$ be a smooth complex projective variety. A $(\Bir(X),X)$-structure on $V$ is a maximal atlas of local charts $\varphi_i:U_i\rightarrow X_i$ such that
\begin{itemize}
	\item the $U_i$ are open sets of $V$ and form a covering;
	\item the $X_i$ are smooth projective varieties birational to $X$;
	\item the $\varphi_i$ are biholomorphic onto their images;
	\item the changes of coordinates $\varphi_i \circ \varphi_j^{-1}:\varphi_j(U_i\cap U_j)\rightarrow\varphi_i(U_i\cap U_j)$ are holomorphic diffeomorphisms which extend to birational maps from $X_j$ to $X_i$. 
\end{itemize}
\end{defi}

\begin{defi}
Let $V$ be a complex manifold. Let $X$ be a smooth complex projective variety. A strict $(\Bir(X),X)$-structure on $V$ is a maximal atlas of local charts $\varphi_i:U_i\rightarrow X$ such that
\begin{itemize}
	\item the $U_i$ are open sets of $V$ and form a covering;
	\item the $\varphi_i$ are biholomorphic onto their images;
	\item the changes of coordinates $\varphi_i \circ \varphi_j^{-1}:\varphi_j(U_i\cap U_j)\rightarrow\varphi_i(U_i\cap U_j)$ are holomorphic diffeomorphisms which extend to birational transformations of $X$. 
\end{itemize}
\end{defi}

\begin{rmk}
Let $X'$ be a smooth birational model of $X$. 
It follows directly from the definition that a $(\Bir(X),X)$-structure on $V$ is the same thing as a $(\Bir(X'),X')$-structure on $V$, and that a strict $(\Bir(X),X)$-structure induces a $(\Bir(X),X)$-structure. But in general it is not true that a strict $(\Bir(X),X)$-structure on $V$ gives rise to a strict $(\Bir(X'),X')$-structure on $V$, see Example \ref{twonotionsaredifferent}.
\end{rmk}

\paragraph{Holonomy and developing map.}
For a classical $(G,X)$-manifold $V$, there exist a group homomorphism $\Hol:\pi_1(V)\rightarrow G$ and a local diffeomorphism $\Dev$ from $\tilde{M}$, the universal cover of $V$, to $X$ such that 
\[\forall \gamma\in\pi_1(V), \Dev\circ \gamma=\Hol(\gamma)\circ \Dev.\]
The map $\Dev$ is called \emph{the developing map} and $\Hol$ is called the \emph{holonomy representation}. A $(G,X)$-structure is uniquely determined by its holonomy and its developing map, up to composition by an element of $G$.

The same proof as in the classical case shows:
\begin{prop}\label{holdev}
Let $X$ be a smooth complex projective variety. Let $V$ be a $(\Bir(X),X)$-manifold. Denote by $\tilde{V}$ the universal cover of $V$ and $\pi$ the quotient map. Fix a base point $v\in V$ and choose a point $w\in\tilde{V}$ such that $\pi(w)=v$. There exist a smooth birational model $Y$ of $X$, a homomorphism $\Hol:\pi_1(V,v)\rightarrow \Bir(Y)$ and a $\pi_1(V,v)$-equivariant meromorphic map $\Dev:\tilde{V}\dashrightarrow Y$ such that 
\[\forall f\in\pi_1(V,v), \Dev\circ f=\Hol(f)\circ \Dev.\]
If $(Y',\Hol',\Dev')$ is another such triple, then there exists a birational map $\sigma$ from $Y$ to $Y'$ such that $\Hol'=\sigma \Hol \sigma^{-1}$ and $\Dev'=\sigma\circ \Dev$. We can choose $(Y,\Hol,\Dev)$ so that $\Dev$ is holomorphic at $w$.
\end{prop}
\begin{pro}
Let $c:[0,1]\rightarrow \tilde{V}$ be a smooth path from $w=c(0)$ to a point $z=c(1)$. The image $c([0,1])$ can be covered by local charts of birational structure $(U_0,\varphi_0:U_0\rightarrow X_0),\cdots,(U_k,\varphi_k:U_k\rightarrow X_k)$ which are pulled-back from local charts on $V$, such that $U_i\cap U_j$ is connected and is non-empty if $j=i+1$. We denote by $g_i$ the map $\varphi_{i-1}\circ \varphi_{i}^{-1}\in\Bir(X_i,X_{i-1})$ which is the unique map such that $g_i\circ \varphi_i$ and $\varphi_{i-1}$ agree on $U_i\cap U_j$; the uniqueness is because of the fact that two birational maps which coincide on a non empty Euclidean open set must be the same. We define $\Dev(z)$ as 
\[\Dev(z)=g_1g_2\cdots g_k\varphi_k(z).\]
To be rigorous, this expression does not associate a value to any point $z$: the $g_i$ are birational so we get only a meromorphic expression. Let us see that $\Dev$ is a well-defined meromorphic map from $\tilde{V}$ to $X_0$; it has milder properties than an arbitrary meromorphic map because locally analytically it behaves as a birational map. The unicity of the $g_i$ guarantees that, once $c$ is fixed, $\Dev$ does not depend on $U_1,\cdots,U_k$, but only on the initial chart $U_0$ at the base point $w$. Choose another path $c'$ from $w$ to $z$. Since $\tilde{V}$ is simply connected, there exists a homotopy $H:[0,1]\times[0,1] \rightarrow \tilde{V}$ between $c$ and $c'$. We can cover $c([0,1]\times[0,1])$ by local charts of birational structure. The uniqueness of the transition maps then shows that $\Dev$ depends only on the homotopy class of $c$.  Around the point $w$, the map $\Dev$ coincides with a coordinate chart, thus is holomorphic. 

Let $f\in\pi_1(V,v)$ be a deck transformation. Let $z=f(w)$ in the above construction of $\Dev$. We can suppose that $U_k=f(U_0)$ and $\varphi_k=\varphi_0\circ f^{-1}$. Then $\Dev(f(w))=g_1g_2\cdots g_k\varphi_0\circ f^{-1}$. Put $\Hol(f)=g_1g_2\cdots g_k$. It belongs to $\Bir(X_0)$. We have $\Dev\circ f=\Hol(f)\circ \Dev$ in a neighbourhood of $w$. Thus $\Dev\circ f=\Hol(f)\circ \Dev$ by analytic continuation.

Let $(Y',\Hol',\Dev')$ be another such triple. Since the set of points of $\tilde{V}$ where a developing map is not defined or is not locally biholomorphic is locally closed of codimension at least one, there exists an open set $U$ of $\tilde{V}$ restricted to which both $\Dev$ and $\Dev'$ are biholomorphic. Then $\Dev\vert_U$ and $\Dev'\vert_U$ are both local birational charts. They have to be compatible, i.e.\ $\Dev'\vert_U\circ (\Dev\vert_U)^{-1}$ extends to a birational map $\sigma$ from $Y$ to $Y'$. By analytic continuation we see that $\sigma$ satisfies $\Hol'=\sigma \Hol \sigma^{-1}$ and $\Dev'=\sigma\circ \Dev$.
\end{pro}

\begin{rmk}
A developing map is locally birational; this means that locally it has a birational expression when written in some complex analytic coordinates. Thus a developing map has no ramification. In particular a ramified covering map is never a developing map.
\end{rmk}

If $V$ is a $(\Bir(X),X)$-manifold, then any finite unramified cover $V'$ of $V$ is equipped with an induced $(\Bir(X),X)$-structure. If $(Y,\Hol,\Dev)$ is a holonomy-developing-map triple for $V$, then the compositions $\pi_1(V')\rightarrow \pi_1(V)\xrightarrow{\Hol}\Bir(Y)$ and $V'\rightarrow V\xrightarrow{\Dev}Y$ form a pair of holonomy and developing map for $V'$.

\begin{prop}\label{devhol}
Let $V$ be a complex manifold with two $(\Bir(X),X)$-structures. Let $(X_1,\Hol_1,\Dev_1)$ and $(X_2,\Hol_2,\Dev_2)$ be pairs of holonomy and developing map associated with these two $(\Bir(X),X)$-structures. The two $(\Bir(X),X)$-structures are the same if and only if there exists $\sigma\in\Bir(X_1,X_2)$ such that $\Hol_2=\sigma \Hol_1 \sigma^{-1}$ and $\Dev_2=\sigma\circ \Dev_1$.
\end{prop}
\begin{pro}
We need to prove the ``if'' part. Let $z$ be a point of the universal cover $\tilde{V}$. Without loss of generality, using the "only if" part (Proposition \ref{holdev}), we can suppose that $\Dev_1$ and $\Dev_2$ are both locally biholomorphic at $z$. Thus on a neighbourhood of $z$, the restrictions of $\Dev_1$ and $\Dev_2$ give local charts for their corresponding birational structures. The hypothesis implies that these charts are compatible, i.e.\ contained in a same maximal atlas. The conclusion follows. 
\end{pro}

\begin{rmk}
Propositions \ref{holdev} and \ref{devhol} hold for strict $(\Bir(X),X)$-structures too. Note that for a strict $(\Bir(X),X)$-structure, the target of the developing map is $X$ itself.
\end{rmk}

The next proposition says that we could alternatively define a birational structure using holonomy and developing map.
\begin{prop}\label{switch}
Let $X$ be a smooth projective variety. Let $V$ be a compact complex manifold and $\tilde{V}$ its universal cover. Let $D:\tilde{V}\dashrightarrow X$ be a meromorphic map that satisfies the following: for every point $w\in\tilde{V}$, there is a Euclidean neighbourhood $W$ of $w$ and a holomorphic diffeomorphism $\varphi$ from $W$ to a Euclidean open set of a birational model $X_w$  of $X$ depending on $w$ such that $D\vert_W\circ \varphi^{-1}$ is the restriction of a birational map. Let $H:\pi_1(V)\rightarrow \Bir(X)$ be a homomorphism of groups such that for every $\gamma\in\pi_1(V)$ we have $H(\gamma)\circ D=D\circ \gamma$. Then $V$ has a $(\Bir(X),X)$-structure for which $(H,D)$ is a holonomy/developing pair. If $X_w=X_0$ are the same for all $w$ then we have a strict $(\Bir(X_0),X_0)$-structure.
\end{prop}
\begin{pro}
Let $v$ be a point of $V$. Choose a point $w\in\tilde{V}$ which projects onto $v$, and a sufficiently small neighbourhood $W$ of $w$ which maps bijectively to a neighbourhood $U$ of $v$. By hypothesis $U$ is biholomorphic to an open subset of a birational model $X_w$. The hypothesis on the local birational property of $D$ and the equivariance of $H$ imply that different choices of $w$ give the same $(\Bir(X),X)$-structure on $U$. Thus $V$ is equipped with a $(\Bir(X),X)$-structure. We leave the reader to verify that $(H,D)$ is indeed a corresponding pair of holonomy and developing map. 
\end{pro}

\begin{cor}\label{switchcor1}
Let $f:X_1\dashrightarrow X_2$ be a birational map between two smooth projective varieties $X_1,X_2$. Let $(X_1,\Hol,\Dev)$ be a holonomy/developing triple associated with a $(\Bir(X_1),X_1)$-structure on a compact complex manifold $V$. Then $(X_2,f\Hol f^{-1},f\circ\Dev)$ is a holonomy/developing triple associated with the same $(\Bir(X_1),X_1)$-structure.
\end{cor}
\begin{pro}
The pair $(f\Hol f^{-1},f\circ\Dev)$ satisfies the conditions of Proposition \ref{switch}, thus defines a $(\Bir(X_2),X_2)$-structure. This $(\Bir(X_2),X_2)$-structure coincides with the original $(\Bir(X_1),X_1)$-structure by Proposition \ref{devhol}.
\end{pro}

\paragraph{Strict $(\Bir(\PP^n),\PP^n)$-structures.}
In general a $(\Bir(X),X)$-structure does not induce a strict $(\Bir(X),X)$-structure. A smooth projective variety birational to $X$ always admits a $(\Bir(X),X)$-structure, but not necessarily a strict $(\Bir(X),X)$-structure as the following example shows:
\begin{exa}\label{twonotionsaredifferent}
Let $X$ be a projective K3 surface. Let $Z$ be the blow-up of $X$ at some point. Suppose by contradiction that $Z$ admits a strict $(\Bir(X),X)$-structure. Consider the developing map $\Dev:Z\dashrightarrow X$. Being locally birational, it induces an injection of the function field of $X$ into that of $Z$. This means that $\Dev$ is a rational dominant map. As a developing map $\Dev$ has no ramification and a K3 surface is simply connected, the degree of $\Dev$ must be one, i.e.\ it is birational. Then we infer that $\Dev$ is the blow-down of the exceptional curve. However by Proposition \ref{holdev} we could choose $\Dev$ so that $\Dev$ is locally biholomorphic around a point on the execptional curve, contradiction.
\end{exa}

The reasoning in the previous example shows:
\begin{lem}\label{devbirational}
Let $X$ be a simply connected smooth variety. For a $(\Bir(X),X)$-structure on $X$, any developing map is birational. Hence the natural $(\Bir(X),X)$-structure on $X$ is the unique one. If $X$ has a strict $(\Bir(Y),Y)$-structure for some $Y$ birational to $X$ then it is the unique strict $(\Bir(Y),Y)$-structure on $X$.
\end{lem}

As rational varieties have the most complicated birational transformation groups, it is natural to ask
\begin{question}\label{strictornot}
\begin{enumerate}
	\item Does a $(\Bir(\PP^n),\PP^n)$-structure always induce a strict $(\Bir(\PP^n),\PP^n)$-structure?
	\item Does every smooth rational variety $X$ of dimension $n$ admit a strict $(\Bir(\PP^n),\PP^n)$-structure?
\end{enumerate}
\end{question}

\begin{prop}\label{itsufficestoverifyavarietyhasastrictstructure}
Questions \ref{strictornot}.1 and \ref{strictornot}.2 are equivalent.
\end{prop}
\begin{pro}
1) implies 2) because every smooth rational variety admits trivially a non-strict $(\Bir(\PP^n),\PP^n)$-structure. Suppose that 2) is true. Let $V$ be a complex manifold with a $(\Bir(\PP^n),\PP^n)$-structure. Let $(U_i,\phi_i:U_i\rightarrow X_i)$ be an atlas for the $(\Bir(\PP^n),\PP^n)$-structure. By hypothesis each $X_i$ has a strict $(\Bir(\PP^n),\PP^n)$-structure. Via $\phi_i$ this equips $U_i$ with a strict $(\Bir(\PP^n),\PP^n)$-structure. Cover $U_i$ by charts $\{U_{ij}\}_j$ of strict $(\Bir(\PP^n),\PP^n)$-structure induced by the one on $X_i$; the $U_{ij}$ are identified via $\phi$ with charts of $(\Bir(\PP^n),\PP^n)$-structure on $X_i$. By Lemma \ref{devbirational} the changes from $U_{ij}$ to $U_{i'j'}$ are birational. Therefore the strict $(\Bir(\PP^n),\PP^n)$-structures on $U_i$ patch together to give a strict $(\Bir(\PP^n),\PP^n)$-structure on $V$ for which the $U_{ij}$ form an atlas. 
\end{pro}

A smooth rational variety $X$ of dimension $n$ is called \emph{uniformly rational} if any $x\in X$ has a Zariski neighbourhood which is isomorphic to a Zariski open set in the affine space $\mathbb{A}^n$. Being rational $X$ has to have such a point; the issue is whether it holds for all points, hence the terminology ``uniformly rational''. Gromov asked:
\begin{question}[Gromov \cite{Gromov89} page 885, see also \cite{BB14}]\label{gromovquestion}
Is every smooth rational complex variety uniformly rational?
\end{question}

It turns out that Gromov's question is equivalent to Question \ref{strictornot} of which the formulation seems not quite algebraic at first glance:
\begin{prop}\label{gromovandbirstructure}
A rational variety of dimension $n$ is uniformly rational if and only if it admits a strict $(\Bir(\PP^n),\PP^n)$-structure.
\end{prop}
\begin{pro}
Suppose that $X$ is a uniformly rational variety of dimension $n$. For any $x\in X$, let $U_x$ be a Zariski neighbourhood isomorphic to an open set of $\mathbb{A}^n$. Then the open sets $U_x$ give an atlas of strict $(\Bir(\PP^n),\PP^n)$-structure.

Now suppose that $X$ admits a strict $(\Bir(\PP^n),\PP^n)$-structure. Let $x\in X$. Take a developing map $\Dev:X\dashrightarrow \PP^n$ which is holomorphic at $x$ (cf. Proposition \ref{holdev}). By Lemma \ref{devbirational} $\Dev$ is birational. Thus $\Dev$ realizes an isomorphism from some Zariski neighbourhood of $x$ to a Zariski open set of $\mathbb{A}^n$.
\end{pro}

Question \ref{gromovquestion} is easy in dimension $1$ or $2$, and is still open in dimension $\geq 3$ (cf. \cite{BB14}). The one dimensional case is trivial because $\PP^1$ is the only smooth rational curve. For completeness we include a proof for the two dimensional case:

\begin{prop}\label{rationalsurfacehasbirationalstructure}
Let $X$ be a smooth rational surface. Then $X$ is uniformly rational and admits a unique strict $(\CRC,\PP^2)$-structure.
\end{prop}
\begin{pro}
Once we prove that $X$ is uniformly rational, we obtain the existence of strict $(\CRC,\PP^2)$-structure by Proposition \ref{gromovandbirstructure} and the uniqueness by Lemma \ref{devbirational}.

A Hirzebruch surface is a $\PP^1$-bundle over $\PP^1$. Cover the base and fiber $\PP^1$ respectively by two pieces of $\mathbb{A}^1$, we see that every Hirzebruch surface can be otbtained by patching four pieces of $\mathbb{A}^2$. Every rational surface different from $\PP^2$ can be obtained from a Hirzebruch surface by blow-ups. The blow-up of $\mathbb{A}^2$ at one point is the union of two Zariski open sets isomorphic to $\mathbb{A}^2$. Hence the uniform rationality.
\end{pro}

So it does no harm if we do not distinguish $(\Bir(\PP^2),\PP^2)$-structure from strict $(\Bir(\PP^2),\PP^2)$-structure. More precisely we have:
\begin{cor}\label{switchcor2}
Let $V$ be a compact complex surface equipped with a $(\Bir(\PP^2),\PP^2)$-structure. Then for any rational surface $X$, there exists a unique strict $(\Bir(X),X)$-structure on $V$ such that, $(\Hol,\Dev)$ is a holonomy/developing pair associated with the strict $(\Bir(X),X)$-structure if and only if $(X,\Hol,\Dev)$ is a holonomy/developing triple associated with the $(\Bir(\PP^2),\PP^2)$-structure.
\end{cor}
\begin{pro}
By Propositions \ref{itsufficestoverifyavarietyhasastrictstructure} and \ref{rationalsurfacehasbirationalstructure} $V$ has a strict $(\Bir(\PP^2),\PP^2)$-structure. Let $(U_i,\phi_i:U_i\hookrightarrow \PP^2)$ be an atlas. The strict $(\Bir(X),X)$-structure on $V$ is constructed as follows. Firstly $\Aut(\PP^2)$ is transitive so any point of $\PP^2$ admits a Zariski neighbourhood which is isomorphic to a Zariski open set of $X$. By considering the intersection of the Euclidean open set $\phi_i(U_i)$ with these Zariski neighbourhoods, we can further subdivide the atlas $(U_i)$ into an atlas $(U_{ij})$ such that there are embeddings $\varphi_{ij}:U_{ij}\hookrightarrow X$ such that $\phi_i\circ\varphi_{ij}^{-1}$ extend to birational maps from $X$ to $\PP^2$. The atlas $(U_{ij}, \varphi_{ij}:U_{ij}\rightarrow X)$ defines the desired strict  $(\Bir(X),X)$-structure.
\end{pro}

\paragraph{Local structure of the developing map.} 
Though the developing map $\Dev$ is in general not holomorphic, it is by construction locally birational. Thus, at least locally, algebro-geometric reasonings could be applied. In dimension two, the indeterminacy set of $\Dev$ is a discrete set of points. We can speak about contracted curves, they are complex analytic subsets of pure dimension $1$. A contracted curve has locally a finite number of components. An irreducible contracted curve is a minimal closed connected $1$-dimensional analytic subset contracted by $\Dev$.

\subsection{Entire curves and Ahlfors-Nevanlinna currents}\label{AhlforsNevanlinnasection}
In this section we give a treatment of families of Ahlfors-Nevanlinna currents which are, we believe, of independent interest. Let $X$ be a smooth projective surface. An \emph{entire curve} on $X$ is a non-constant holomorphic map $\xi:\CC\rightarrow X$. An entire curve $\xi$ is called \emph{transcendental} if its image is not contained in an algebraic curve of $X$. We can associate to a transcendental entire curve $\xi$ a (a priori non unique) closed positive current, called \emph{Ahlfors-Nevanlinna current}. We need to prove a variant of the construction for a family of entire curves. We recall first the process for a single entire curve (see \cite{Mcq98}, \cite{Bru99} for Ahlfors-Nevanlinna currents and \cite{Dem97} for the functions we use below).

For a differential form $\eta\in \mathcal{A}^2(X)$ and for $r>0$ we put 
\[T_{\xi,r}(\eta)=\int_0^r\frac{ds}{s}\int_{\DD_s}\xi^*\eta\]
where $\DD_s\subset\CC$ is the disk of radius $s$. We fix a \Kah{} form $\omega\in \mathcal{A}^{1,1}(X)$. Consider the positive currents defined by 
\[
\Phi_r(\eta)=\frac{T_{\xi,r}(\eta)}{T_{\xi,r}(\omega)}, \quad \forall \eta\in \mathcal{A}^2(X).
\]
The family $\{\Phi_r\}_{r>0}$ is bounded, so we can find a sequence of radii $(r_n)_{n\in\NN}$ such that $r_n\rightarrow +\infty$ and $\Phi_{r_n}$ converges weakly to a positive current $\Phi$. For the limit $\Phi$ to be a closed current, we need a smart choice of the sequence $(r_n)$. Let us denote by $A(r)$ the area of $\xi(\DD_r)$ and $L(r)$ the length of $\xi(\partial \DD_r)$ with respect to the Riemannian metric induced by $\omega$. Then $T_{\xi,r}(\omega)$ may be written as 
\[T_{\xi,r}(\omega)=\int_0^r A(s)\frac{ds}{s}.\]
We have 
\[
\operatorname{limsup}_{r\rightarrow \infty}\frac{T_{\xi,r}(\omega)}{\log r}=\infty
\]
since $\xi$ is transcendental (cf. \cite{Dem97}).
We define
\[S_{\xi,r}(\omega)=\int_0^r L(s)\frac{ds}{s}.\]
For $\beta\in \mathcal{A}^1(X)$, Stokes' theorem and the compactness of $X$ imply the inequality
\[
\vert T_{\xi,r}(d \beta)\vert \leq \int_0^r\frac{ds}{s}\int_{\partial \DD_s} \vert \xi^*\beta \vert \leq constant\cdot S_{\xi,r}(\omega),
\]
where the constant on the right side depends on $\beta$ but not on $r$. Therefore to obtain a closed limit current $\Phi$, we need a sequence of radii $(r_n)_n$ such that 
\[
\frac{S_{\xi,r_n}(\omega)}{T_{\xi,r_n}(\omega)}\rightarrow 0, \quad \text{when}\ n\rightarrow \infty.
\]
The existence of such a sequence of radii is guaranteed by the following lemma (see \cite{Bru99}):
\begin{lem}[Ahlfors \cite{Ahl35}]
Let $R>0$, $\epsilon >0$ be two positive real numbers. Denote by $B(\xi,\epsilon)$ the set $\{r>R\vert S_{\xi,r}(\omega)>\epsilon T_{\xi,r}(\omega)\}$. Then
\[\int_{B(\xi,\epsilon)}\frac{dr}{r\log r}<\infty.\]
In particular $\operatorname{liminf}_{r\rightarrow \infty}\frac{S_{\xi,r}(\omega)}{T_{\xi,r}(\omega)}=0$.
\end{lem}
Note that the measure of $(R,\infty)$ with respect to $\frac{dr}{r\log r}$ is infinite, so the above lemma implies that we can choose an appropriate sequence of radii simultaneously for a finite number of entire curves:
\begin{lem}\label{afinitenumberofentirecurves}
Let $\xi_1,\cdots,\xi_k$ be $k$ transcendental entire curves on $X$. There exists a sequence of radii $(r_n)_{n\in\NN}$ such that for each $i\in\{1,\cdots,k\}$, the sequence $\left(\frac{T_{\xi_i,r_n}(\cdot)}{T_{\xi_i,r_n}(\omega)}\right)_{n\in\NN}$ converges weakly to a closed positive current.
\end{lem}

A closed positive current $\Phi$ constructed by the above limit process is called an \emph{Ahlfors-Nevanlinna current} associated with the entire curve $\xi$, it depends on the choice of a sequence of radii $(r_n)_n$. 

A cohomology class is called \emph{nef} if its intersections with all curves are non negative. We refer the reader to \cite{Mcq98}, \cite{Bru99} for the following:
\begin{lem}\label{ANcurrentisnef}
Let $[\Phi]\in \operatorname{H}^{1,1}(X,\RR)$ be the cohomology class of an Ahlfors-Nevanlinna current associated with a transcendental entire curve. Then $[\Phi]$ is nef. In particular $[\Phi]^2\geq 0$.
\end{lem}

We will need to consider some families of entire curves. To treat the Ahlfors-Nevanlinna currents simultaneously in family, we need some control on the variation of entire curves. The following very restricted notion will be sufficient for our proof.
\begin{defi}\label{defuniformfamily}
A family of entire curves parametrized by a real manifold $B$ is a differentiable map $B\times \CC\rightarrow X, (b,z)\mapsto \xi_b(z)$ such that $\xi_b$ is an entire curve for every $b\in B$. A family of entire curves $(\xi_b)_{b\in B}$ is called uniform if the following condition is satisfied: $\forall b_0\in B, \forall \delta>0$, there exists a neighborhood $U$ of $b_0$ such that 
\[
\forall b\in U, \forall z\in \CC, \left\vert \vert\xi_b'(z)\vert-\vert\xi_{b_0}'(z)\vert \right\vert<\delta \vert \xi_{b_0}'(z) \vert
\]
where the absolute values are measured with respect to a fixed \Kah{} metric on $X$. In other words a family of entire curves is uniform if nearby pull-backed metrics are close in proportion.
\end{defi}

There exist non-trivial uniform families of transcendental entire curves on complex projective surfaces, for example there exist families of Levi-flat hypersurfaces foliated by entire curves (see Remark 1.6 of \cite{Der05}). Our interest in this notion is explained by the following lemma:
\begin{lem}\label{lemfamilyofentirecurves}
Let $(\xi_b)_{b\in B}$ be a uniform family of transcendental entire curves on $X$. Let $A$ be a compact $C^{\infty}$-path connected subset of $B$. Then there exists a sequence of radii $(r_n)_{n\in\NN}$ using which an Ahlfors-Nevanlinna current associated with $\xi_a$ can be constructed for all $a\in A$. After fixing such a sequence $(r_n)_{n\in\NN}$, the Ahlfors-Nevanlinna currents associated with the $\xi_a$ all have the same cohomology class. 
\end{lem}
\begin{pro}
We prove first that there exists a common choice of the sequence of radii. Let us fix a sufficiently small real number $\delta>0$. By the definition of uniform family and by compactness of $A$, we can find a finite number of points $a_1,\cdots,a_k$ in $A$ and their neighbourhoods $U_1,\cdots,U_k$ in $B$ such that 
\begin{itemize}
	\item $\forall i\in\{1,\cdots,k\}, \forall a\in U_i, \forall z\in \CC, \left\vert \vert\xi_a'(z)\vert-\vert\xi_{a_i}'(z)\vert \right\vert<\delta \vert \xi_{a_i}'(z) \vert$;
	\item $A\subset \cup U_i$.
\end{itemize}

By lemma \ref{afinitenumberofentirecurves}, we can take a sequence of radii $(r_n)_{n\in\NN}$ that works for all the $\xi_{a_i}, 1\leq i\leq k$. We denote by $\lambda$ the Lebesgue measure on $\CC$. Let $a\in U_i$. We have
\begin{align*}
T_{\xi_a,r}(\omega)=\int_0^r A(s)\frac{ds}{s}=\int_0^r \int_{\DD_r}\vert \xi_a'(z) \vert ^2 d\lambda(z) \frac{ds}{s}
\end{align*}
and            
\begin{align*}
\vert T_{\xi_a,r}(\omega)-T_{\xi_{a_i},r}(\omega) \vert & \leq \int_0^r \int_{\DD_r}\vert \vert \xi_a'(z) \vert ^2 -\vert \xi_{a_i}'(z) \vert ^2 \vert d\lambda(z) \frac{ds}{s}\\
& \leq (2\delta +\delta^2)\int_0^r \int_{\DD_r}\vert \xi_{a_i}'(z) \vert ^2 d\lambda(z) \frac{ds}{s}\\
& =(2\delta +\delta^2) T_{\xi_{a_i},r}(\omega) 
\end{align*}
Similarly we have
\begin{align*}
\vert S_{\xi_a,r}(\omega)-S_{\xi_{a_i},r}(\omega) \vert \leq \delta S_{\xi_{a_i},r}(\omega).
\end{align*}
Consequently 
\[
\frac{S_{\xi_a,r}(\omega)}{T_{\xi_a,r}(\omega)}\leq \frac{1+\delta}{1-2\delta-\delta^2}\frac{S_{\xi_{a_i},r}(\omega)}{T_{\xi_{a_i},r}(\omega)}.
\]
In particular we have 
\[
\lim_{n\rightarrow \infty}\frac{S_{\xi_a,r_n}(\omega)}{T_{\xi_a,r_n}(\omega)}=\lim_{n\rightarrow \infty}\frac{S_{\xi_{a_i},r_n}(\omega)}{T_{\xi_{a_i},r_n}(\omega)}=0
\]
so that the sequence $(r_n)_n$ can be used to construct Ahlfors-Nevanlinna currents for all $a\in A$. Hence we can talk about \emph{the} Ahlfors-Nevanlinna currents $\Phi_a$ associated with the $\xi_a$ and this fixed sequence $(r_n)_n$. 

Let $a,b$ be two points in $A$. We now prove that the Ahlfors-Nevanlinna currents $\Phi_a,\Phi_b$ are cohomologous. It is sufficient to treat the case where $a=a_i$ and $b\in U_i$. Take a $C^{\infty}$-path $c:[0,1]\rightarrow U_i$ such that $c(0)=a,c(1)=b$. We denote by $F$ the induced map $[0,1]\times \CC\rightarrow X,F(s,z)=\xi_{c(s)}(z)$. Let $\eta\in A^2(X)$. Applying Stokes' theorem, we have
\begin{equation}
T_{\xi_a,r}(\eta)-T_{\xi_b,r}(\eta)=\int_0^r\int_{[0,1]\times \DD_s}F^*(d\eta)\frac{ds}{s}-\int_0^r\int_{[0,1]\times \partial\DD_s}F^*(\eta)\frac{ds}{s}. \label{eq:stokesformula}
\end{equation}
We denote by $\Theta_r$ the current of dimension $3$ defined by $\Theta_r(\beta)=\int_0^r\int_{[0,1]\times \DD_s}F^*(\beta)\frac{ds}{s}$ for $\beta\in A^3(X)$, and by $\Xi_r$ the current of dimension $2$ defined by $\Xi_r(\beta)=\int_0^r\int_{[0,1]\times \partial\DD_s}F^*(\beta)\frac{ds}{s}$. We have
\[
\frac{T_{\xi_a,r}(\eta)}{T_{\xi_a,r}(\omega)}-\frac{T_{\xi_b,r}(\eta)}{T_{\xi_b,r}(\omega)}
=\frac{T_{\xi_a,r}(\eta)-T_{\xi_b,r}(\eta)}{T_{\xi_b,r}(\omega)}+\frac{T_{\xi_b,r}(\omega)-T_{\xi_a,r}(\omega)}{T_{\xi_a,r}(\omega)T_{\xi_b,r}(\omega)}T_{\xi_b,r}(\eta)
\]
which with Equation \eqref{eq:stokesformula} implies
\begin{equation}
\frac{T_{\xi_a,r}(\eta)}{T_{\xi_a,r}(\omega)}-\frac{T_{\xi_b,r}(\eta)}{T_{\xi_b,r}(\omega)}
=\frac{1}{T_{\xi_b,r}(\omega)}\left( d\Theta_r(\eta)-\Xi_r(\eta) \right)+\frac{T_{\xi_b,r}(\omega)-T_{\xi_a,r}(\omega)}{T_{\xi_a,r}(\omega)}\frac{T_{\xi_b,r}(\eta)}{T_{\xi_b,r}(\omega)}.\label{eq:bigformula}
\end{equation}

We want to show that along the sequence of radii $(r_n)_{n\in \NN}$, the right side of Equation \eqref{eq:bigformula} converges weakly to an exact current. We first estimate $\Xi_r$. By compactness of $X$, we have 
\[
\vert \Xi_r(\eta)\vert=\left\vert \int_0^r\int_{[0,1]\times \partial\DD_s}F^*(\eta)\frac{ds}{s} \right\vert \leq M(\eta) \int_0^r\int_{[0,1]\times \partial\DD_s}\vert F^*(\omega)\vert\frac{ds}{s}
\]
where $M(\eta)$ is a constant that depends on $\eta$ but not on $r$. By Fubini's theorem, we deduce from the above inequality that
\begin{align*}
\vert \Xi_r(\eta)\vert & \leq M(\eta) \int_0^r\int_0^1 L(c(u),s) du \frac{ds}{s}
\end{align*}
where $L(c(u),s)$ is the length of $\xi_{c(u)}(\partial \DD_s)$ with respect to the \Kah{} metric defined by $\omega$. Using the fact that the path $c$ lies in $U_i$, we have further
\begin{align*}
\vert \Xi_r(\eta)\vert & \leq M(\eta) \int_0^r\int_0^1 (1+\delta)L(a,s) du \frac{ds}{s}=M(\eta)(1+\delta)S_{\xi_a,r}(\omega). 
\end{align*}
This implies that the sequence of currents $(\Xi_{r_n}/T_{\xi_b,r_n}(\omega))_{n\in \NN}$ converges weakly to $0$.

Now we estimate the last term of Equation \eqref{eq:bigformula}. By Stokes' formula, we have
\begin{equation}
T_{\xi_b,r}(\omega)-T_{\xi_a,r}(\omega)=\int_0^r\int_{[0,1]\times \DD_s}F^*(d\omega)\frac{ds}{s}-\int_0^r\int_{[0,1]\times \partial\DD_s}F^*(\omega)\frac{ds}{s}. \label{eq:kahlerformisclosed}
\end{equation}
Since the form $\omega$ is \Kah{}, we have $d\omega=0$ and the first term of the right side of \eqref{eq:kahlerformisclosed} vanishes. The second term at the right side of \eqref{eq:kahlerformisclosed} is dominated by $S_{\xi_a,r}(\omega)$. It follows immediately that the last term of \eqref{eq:bigformula} converges weakly to zero along the sequence of radii $(r_n)_n$. 

Finally we estimate $\Theta_r$. Note that, since the other terms in Equation \eqref{eq:bigformula} all converge weakly along the sequence $(r_n)_n$, the sequence $(d\Theta_{r_n}/T_{\xi_b,r_n}(\omega))_n$ converges weakly too. However this does not imply that $(\Theta_{r_n}/T_{\xi_b,r_n}(\omega))_n$ converges weakly. Again using Fubini's theorem and compactness of $X$, we have
\[
\vert \Theta_r(\beta) \vert\leq N(\beta) T_{\xi_b,r}(\omega)
\]
where $N(\beta)$ is a constant which depends on $\vert\beta\vert$ and on $\delta$ but not on $r$. Thus the $\Theta_{r_n}/T_{\xi_b,r_n}(\omega)$ form a bounded family and there exists a subsequence $(r_{n_j})_j$ of $(r_n)_n$ such that $\Theta_{r_{n_j}}/T_{\xi_b,r_{n_j}}(\omega)$ converges weakly to a current $\Theta$. Hence, the weak limit of $(d\Theta_{r_n}/T_{\xi_b,r_n}(\omega))_n$ is exact because
\[
\lim_{n\rightarrow \infty}\frac{d\Theta_{r_n}(\eta)}{T_{\xi_b,r_n}(\omega)}=\lim_{j\rightarrow \infty}d\left(\frac{\Theta_{r_{n_j}}}{T_{\xi_b,r_{n_j}}(\omega)}\right)(\eta)=d\Theta(\eta).
\]
The conlusion follows.
\end{pro}
Note that to construct Ahlfors-Nevanlinna currents it is not necessary for $\omega$ to be \Kah{}: a Hermitian metric would be sufficient. However the property $d\omega=0$ is used in the last part of the above proof.

\subsection{Transverse invariant measures}
All the materials in this section can be found in \cite{Ghy99} and \cite{FS08}. Let $M$ be a compact Hausdorff topological space. A structure of \emph{lamination by Riemann surfaces} on $M$ is an atlas $\mathcal{L}$ of charts $h_i:U_i\rightarrow \DD\times B_i$ where $\DD$ is the unit disk in $\CC$, the $B_i$ are topological spaces, the $h_i$ are homeomorphisms and the $U_i$ are open sets of $M$ which cover $M$; the changes of coordinates $h_{ij}=h_j\circ h_i^{-1}$ are of the form $(f_{ij}(z,b),g_{ij}(b))$ where the $f_{ij}$ are holomorphic in $z$ and the $g_{ij}$ are continuous. A connected component of $V_c=\{(z,b)\vert b=c\}$ in a chart $U_i$ is called a \emph{plaque}. A minimal connected subset of $M$ which contains all plaques that it intersects is called \emph{a leaf}. A lamination by Riemann surfaces $(M,\mathcal{L})$ is \emph{transversally smooth} if the $B_i$ are real manifolds and if the $g_{ij}$ are smooth maps. A \emph{transverse invariant measure} $\mu$ on $(M,\mathcal{L})$ is a family of locally finite positive measures $\mu_i$ on the topological spaces $B_i$ such that if $B\subset B_i$ is a measurable set contained in the domain of definition of $g_{ij}$, then $\mu_i(B)=\mu_j(g_{ij}(B))$. 

From now on we make the hypothesis that $(M,\mathcal{L})$ is a lamination by Riemann surfaces contained in a \Kah{} surface $X$. This hypothesis is just for convenience of the presentation and everything we will say makes sense without assuming that there is an ambient surface $X$. Examples to keep in mind are Levi-flat hypersurfaces and saturated sets of holomorphic foliations. We say that a continuous $(1,0)$-form $\beta$ on $X$ \emph{defines the lamination} $(M,\mathcal{L})$ if $\beta\wedge [V_c]=0$ for every plaque $V_c$, where $[V_c]$ is the current of integration on the plaque $V_c$. A closed positive current $\Theta$ of bidimension $(1,1)$ on $X$ is \emph{directed by $(M,\mathcal{L})$} if it is supported on $M$ and if $\Theta\wedge \beta=0$ for all $\beta$ defining $(M,\mathcal{L})$. Our purpose of introducing the above notions is the following theorem:
\begin{thm}[Sullivan \cite{Sul76}]\label{Sullivanthm}
Let $(M,\mathcal{L})$ be a transversally smooth lamination by Riemann surfaces contained in a \Kah{} surface $X$. A transverse invariant measure on $(M,\mathcal{L})$ is the same thing as a closed positive current directed by $(M,\mathcal{L})$ via the following correspondence: in a chart $h_i:U_i\rightarrow \DD\times B_i$, a closed positive directed current $T$ may be written as
\[T=\int_{B_i} [V_b]d\mu(b)\] where $\mu$ is a transverse invariant measure and the $[V_b]$ are integrations on plaques.
\end{thm}
We will apply Sullivan's theorem to Ahlfors-Nevanlinna currents associated with entire curves tangent to the lamination, thanks to the following construction studied by Plante:
\begin{thm}[Plante \cite{Pla75}, see also \cite{Ghy99}, \cite{FS08}]\label{Plantethm}
Let $(M,\mathcal{L})$ be a lamination by Riemann surfaces contained in a \Kah{} surface $X$. Let $f:\CC\rightarrow X$ be a transcendental entire curve contained in a leaf of the lamination and let $\Phi_f$ be an Ahlfors-Nevanlinna current associated with $f$. Then $\Phi_f$ is directed by $(M,\mathcal{L})$.
\end{thm}

\section{Inoue surfaces of type $S^0$.}\label{secondsection}
\subsection{Description}
Let $M\in \operatorname{SL} _3(\ZZ)$ be a matrix with eigenvalues $\alpha,\beta,\bar{\beta}$ such that $\alpha>1$ and $\beta\neq \bar{\beta}$. Note that $\alpha$ is irrational and $\vert\beta\vert<1$. We choose a real eigenvector $(a_1,a_2,a_3)$ corresponding to $\alpha$ and a complex eigenvector $(b_1,b_2,b_3)$ corresponding to $\beta$. Let $G_M$ be the subgroup of $\Aut(\PP^1\times\PP^1)$ generated by 
\begin{align*}
g_0:&(x,y)\mapsto (\alpha x,\beta y)\\
g_i:&(x,y)\mapsto (x+a_i,y+b_i)\quad \text{for} \ i=1,2,3.
\end{align*}
Denote by $\HH$ the upper half plane, viewed as an open subset of $\PP^1=\CC\cup \{\infty\}$. The action of $G_M$ preserves $\HH\times\CC$; it is free and properly discontinuous. The quotient $S_M=\HH\times\CC/G_M$ is a compact non-\Kah{} surface without curves called an \emph{Inoue surface of type $S^0$} (\cite{Ino74}). Note that we should have included the choices of $(a_1,a_2,a_3)$ and $(b_1,b_2,b_3)$ in the notation of $S_M$. By construction it has an $(\Aff _2(\CC),\CC^2)$-structure where by $\Aff _2(\CC)$ we denote the affine transformation group of $\CC^2$. In particular an Inoue surface of type $S^0$ has a natural $(\CRC,\PP^2)$-structure. 

Consider the following solvable Lie group which is a subgroup of $\Aff _2(\CC)$:
\[
\operatorname{Sol}^0=\left\{\begin{pmatrix}\vert\lambda\vert^{-2} &0&a\\0&\lambda &b\\0&0&1\end{pmatrix},\lambda \in\CC^*,a\in\RR,b\in\CC\right\}.
\]
The group $\operatorname{Sol}^0$ is a semi-direct product $(\CC\times\RR)\rtimes\CC^*$. It acts transitively on $\HH\times\CC$; the stabilizer of a point is isomorphic to $\mathbb{S}^1$. The group $G_M$ defining the Inoue surface $S_M$ is a lattice in $\operatorname{Sol}^0$; conversely any torsion free lattice of $\operatorname{Sol}^0$ gives an Inoue surface of type $S^0$. The three elements $g_1,g_2,g_3$ generate a free abelian group of rank three; denote it by $A_M$. The group $G_M$ is a semi-direct product $A_M\rtimes \ZZ$ where the $\ZZ$ factor is generated by $g_0$. We have $g_0g_ig_0^{-1}=g_1^{m_{i1}}g_2^{m_{i2}}g_3^{m_{i3}}$ where the $m_{ij}$ are the entries of the matrix $M$. Note that a finite unramified cover of an Inoue surface of type $S^0$ is an Inoue surface of type $S^0$. 

The following lemma says that $G_M$ has few normal subgroups. In particular the commutator $[G_M,G_M]$ is a finite index subgroup of $A_M$. 

\begin{lem}\label{kernellem}
If $K$ is a non-trivial normal subgroup of $G_M$, then either $K$ is of finite index in $A_M$ or $K$ is of finite index in $G_M$.
\end{lem}
\begin{pro}
The conjugation action of $g_0$ on $A_M$ is just the action of $M\in\operatorname{SL} _3(\ZZ)$ on $\ZZ^3$. For all $v\in\ZZ^3\backslash\{0\}$, the iterates $M^n v$ generate a finite index subgroup of $\ZZ^3$. Thus, if $K\cap A_0$ is non trivial, then $K\cap A_0$ is a free $\ZZ$-module of rank $3$ and is of finite index in $A_M$. To conclude, we need only remark that, by the semi-direct product structure, the intersection of a normal subgroup of $G_M$ with $A_M$ cannot be trivial.
\end{pro}

\begin{lem}\label{injectionintoPGLtC}
Let $\sigma:G_M\rightarrow \PGL_2(\CC)$ be an injective morphism. Then for some affine coordinate $\PP^1=\{x\in\CC\}\cup \{\infty\}$, the images $\sigma(g_i), i=0,\cdots,3$, viewed as homographies of $\PP^1$, may be written as
\begin{align*}
\sigma(g_i):\ & x\mapsto x+u_i, \quad i=1,2,3\\ 
\sigma(g_0):\ & x\mapsto \nu x
\end{align*}
for some $\nu,u_i\in\CC^*$.
\end{lem}
\begin{pro}
As $\sigma(g_i),i=1,2,3$ commute with each other, we have two possibilities for them: we can find an affine coordinate $x$ such that they are either $x\mapsto x+u_i$ with $u_i\neq 0$ or $x\mapsto \alpha_i x$ with $\alpha_i$ of infinite order. 

Suppose by contradiction that the $\sigma(g_i)(x)=\alpha_i x$. Since $A_M$ is normal, $\sigma(g_0)$ preserves the set of fixed points of $\sigma(A_M)$, which is $\{0,\infty\}$. Hence $\sigma(g_0)(x)=\gamma x^{\pm 1}$. But then the action of $\sigma(g_0)$ on $\sigma(A_M)$ has finite order, a contradiction.

Hence the $\sigma(g_i)$ are $x\mapsto x+u_i$. The invariance of the fixed point $\infty$ implies that $\sigma(g_0)$ is $x\mapsto\nu x+\delta$ where $\nu$ satisfies that $\nu u_i=m_{i1}u_1+m_{i2}u_2+m_{i3}u_3$ and $\delta$ is arbitrary. Then the change of coordinates $x\mapsto x-\frac{\delta}{\nu}$ allows us to write the $\sigma(g_i)$ as in the statement of the lemma.
\end{pro}

\begin{lem}\label{theonlyfoliations}
The only (possibly singular) holomorphic foliations on $S_M$ are the two obvious ones coming from the horizontal foliation and the vertical foliation of $\HH\times \CC$. 
\end{lem}
\begin{pro}
This is already observed by Brunella in \cite{Bru97} without proof details. Here we give a proof for completeness. See \cite{Bru00} for the terminology we use concerning holomorphic foliations. Suppose by contradiction that $\FF$ is a non-necessarily saturated holomorphic foliation on $S_M$ different from the two obvious ones. Since $S_M$ has no curves (\cite{Ino74}), the singularities of $\FF$ are necessarily isolated. We compare $\FF$ with one of the two obvious foliations: the tangency locus is empty because otherwise it would be a curve on $S_M$. Since the tangency locus contains the singularities of $\FF$, we deduce that $\FF$ is a regular holomorphic foliation, transverse to the two obvious foliations. 

We denote by $T$ the tangent bundle of $S_M$, by $T^*$ its dual and by $K$ the canonical bundle of $S_M$. We denote by $F_0$ the normal bundle of one obvious foliation; the normal bundle of the other obvious foliation is then $-K-F_0$ (here we use notations of \cite{Ino74}). Let $F$ be the normal bundle of $\FF$. The foliation $\FF$ corresponds to a non-zero global section of $T^*\otimes F$. It is proved in \cite{Ino74} (see the first two sentences of sections 6 and 8) that the only line bundles $F$ on $S_M$ such that $T^*\otimes F$ has non-zero sections are $F=F_0$ or $F=-F_0-K$. In other words, $\FF$ and one of the obvious foliations share the same normal bundle. As $\FF$ is everywhere transverse to this foliation, the two sections of $T^*\otimes F$ corresponding to the two foliations trivialize the sheaf $T^*\otimes F$, i.e.\ $T^*$ is isomorphic to $(-F)\oplus (-F)$. However $T^*=(-F_0)\oplus (F_0+K)$ and $F$ is either $F_0$ or $-F_0-K$. This leads to a contradiction as $K$ is not trivial.
\end{pro}

The surface $\overline{S_M}=\HH\times\CC/A_M$ is an infinite cyclic cover of $S_M$. As a real manifold, $\overline{S_M}$ admits a fibration $\rho:\overline{S_M}\rightarrow \RR_+^{*}$ where $\RR_+^{*}=\{t\sqrt{-1},t\in\RR_+^*\}$ is the vertical axis of the component $\HH$ of $\HH\times\CC$. The fibers of $\rho$, denoted by $F_t$, are quotients of $\{x+t\sqrt{-1},x\in\RR\}\times \CC$ by $A_M$; they are real tori of dimension $3$. The $F_t$ are Levi-flat hypersurfaces in $\overline{S_M}$ and they are foliated by entire curves coming from the vertical complex lines in $\HH\times\CC$. 

\begin{lem}\label{theonlymeasure}
Up to multiples, there is only one transverse invariant measure on the Levi flat hypersurface $F_t$.
\end{lem}
\begin{pro}
Recall that $(a_1,a_2,a_3)$ is an eigenvector associated with the irrational eigenvalue $\alpha$ of $M\in\operatorname{SL} _3(\ZZ)$. A transverse invariant measure on $F_t$ is induced by a measure on $\RR=\{x+t\sqrt{-1},x\in\RR\}$ which is invariant under the group of translations generated by $x\mapsto x+a_i,i=1,2,3$. This latter group is a dense subgroup of $\RR$ so the transverse invariant measure must be a multiple of the Lebesgue measure. 
\end{pro}

\begin{lem}\label{nomeasures}
The two obvious foliations on $S_M$ are not transversely Euclidean.
\end{lem}
\begin{pro}
The two dimensional Euclidean isometry group is the semi-direct product $\RR^2\rtimes\mathrm{SO}_2(\RR)$, where $\RR^2$ is the group of translations and $\mathrm{SO}_2(\RR)$ is the group of rotations. Suppose by contradiction that one obvious foliation is transversely Euclidean. Then to this transverse Euclidean structure are associated a holonomy representation $\rho:G_M\rightarrow \RR^2\rtimes\mathrm{SO}_2(\RR)$ and a continuous $\rho$-equivariant developing map $D:T\rightarrow \RR^2$, where the space of leaves $T$ is $\HH$ or $\CC$ depending on which of the two obvious foliations we are looking at. We prove first that $\rho$ is injective by contradiction. Suppose that the kernel $K$ of $\rho$ is not trivial, then it is a finite index subgroup of $A_M$ by Lemma \ref{kernellem}. As $A_M$ is a group of translations on $T$ which is isomorphic to $\ZZ^3$, the closure of any $A_M$-orbit contains at least one real line (for $T=\HH$ a subgroup of $A_M$ isomorphic to $\ZZ^2$ would be sufficient as the $a_i$ are real).  The same holds for $K$-orbits. Then by the $\rho$-equivariance and the continuity of the developing map, $D$ is constant on each of these real lines. This contradicts the fact that the developing map is locally homeomorphic. 

We know now that $\rho$ is injective. As $A_M$ is abelian, we must have $\rho(A_M)\subset \RR^2$. The conjugation action of $g_0$ on $\RR^3=A_M\otimes \RR$ and that of $\rho(g_0)$ on $\RR^2$ are linear maps. We think of $g_0$ and $\rho(g_0)$ as linear maps via their conjugation actions. The group morphism $\rho$ induces a linear map $\pi:\RR^3\rightarrow \RR^2$ which is equivariant under the actions of $g_0$ and $\rho(g_0)$, i.e.\ we have $\pi\circ g_0=\rho(g_0)\circ \pi$. This is not possible because $\rho(g_0)$ is a rotation while $g_0$ corresponds to the matrix $M$ whose eigenvalues are $\alpha,\beta,\overline{\beta}$ with $\alpha>1$ and $\vert \beta \vert<1$.
\end{pro}

\subsection{Proof of Theorem \ref{mainthm} for Inoue surfaces of type $S^0$}\label{proof1}
Let $S_M$ be an Inoue surface of type $S^0$. We fix a $(\Bir(X),X)$-structure on $S_M$ where $X$ is some projective surface. We want to prove that $X$ is rational and the structure is just the obvious affine structure. Let $(Y,\Dev,\Hol)$ be a corresponding holonomy/developing triple as in Proposition \ref{holdev}. We will denote by $\pi$ the covering map from $\HH\times\CC$ to $S_M$.

Lemma \ref{kernellem} says that there are only three possibilities for the holonomy representation. It is easy to rule out the first possibility: if the holonomy had finite image then the developing map would induce a meromorphic locally birational map from a finite unramified cover of $S_M$ to $Y$, contradicting the fact that $S_M$ has algebraic dimension zero. The second possibility is that the kernel $K$ of the holonomy is a finite index subgroup of $A_M$. Then $K\rtimes \ZZ$ has finite index in $G_M$; in this case by considering the corresponding finite unramified cover of $S_M$ and the induced birational structure, we can suppose that $K=A_M$. We will prove in a first step that this case is not possible either. Then we examine the last possibility where the holonomy representation is injective.

\subsubsection{The holonomy is not cyclic}\label{subsectionnoncylicholonomy}
The proof of the following proposition will occupy the rest of this section.
\begin{prop}\label{propnoncyclic}
The image of the holonomy representation is not cyclic.
\end{prop}
We want to prove it by contradiction. \emph{We can and will assume in the sequel that the kernel of $\Hol$ is exactly $A_M$.} Thus the developing map $\Dev:\HH\times\CC\dashrightarrow Y$ factorizes through $\overline{\Dev}:\overline{S_M}\dashrightarrow Y$. We will call the latter map the developing map too.
\begin{lem}\label{finitelymanycontractedcurves}
The developing map $\overline{\Dev}$ has only a finite number of irreducible contracted curves.
\end{lem}
\begin{pro}
Consider the (real) fibration $\rho:\overline{S_M}\rightarrow \RR^{*+}$. The fibers $F_t$ are compact and $\overline{dev}$ is locally birational, so each fiber intersects only a finite number of irreducible contracted curves. Thus it is sufficient to prove that every irreducible contracted curve intersects all the fibers. In other words, let $C\subset \overline{S_M}$ be an irreducible contracted curve, then we want to prove that $\rho(C)=\RR^{*+}$. Since $\rho$ is proper and $C$ is closed, the image $\rho(C)$ is closed in $\RR^{*+}$. It is then sufficient to prove that $\rho(C)$ is open. For this purpose it is more convenient to look at the universal covering $\HH\times\CC$. Let $\tilde{C}\subset \HH\times\CC$ be a component of the inverse image of $C$. Then the projection of $\tilde{C}$ onto $\HH$ is not a point because $C$ cannot be contained in a leaf of the foliation. Thus the projection is open since a holomorphic map is open. Therefore $\rho(C)$, identified as the projection of $\tilde{C}$ onto the vertical axis of $\HH$, is also open.
\end{pro}

Let $C\subset \overline{S_M}$ be an irreducible contracted curve and $q\in Y$ be the point onto which $C$ is contracted. Take a point $c\in C$ which is not an indeterminacy point of $\overline{\Dev}$. Take a chart of birational structure $U\subset \overline{S_M}$ at $c$ so that the restriction $\overline{\Dev}\vert_U$ is analytically equivalent to a birational map. By Zariski's decomposition of birational maps, we can blow up $Y$ at $q$ and its infinitely near points to obtain a surface $Y'$ such that the map $U\dashrightarrow Y'$ induced by $\overline{\Dev}\vert_U$ does not contract $C\cap U$. By analytic continuation, the map $\overline{S_M}\dashrightarrow Y'$ induced by $\overline{\Dev}$ does not contract $C$. As $\overline{\Dev}$ has only finitely many irreducible contracted curves by Lemma \ref{finitelymanycontractedcurves}, by repeating the above process we can find a rational surface $Y^*$, obtained by blowing up $Y$ a finite number of times, such that the induced map $\overline{S_M}\dashrightarrow Y^*$ has no contracted curves. By replacing $Y$ with $Y^*$ \emph{we will suppose from now on that $\Dev$ and $\overline{\Dev}$ have no contracted curves}. 

However $\overline{\Dev}$ may still have indeterminacy points. We denote by $I\subset \overline{S_M}$ the indeterminacy set of $\overline{\Dev}$, it is a discrete set. The map $\overline{\Dev}:\overline{S_M}\dashrightarrow Y$ is locally biholomorphic outside $I$. We will call $\overline{\Dev}(\overline{S_M}\backslash I)$ the \emph{image} of $\overline{S_M}$ (it is also the image of $\Dev$). 

The deck transformation group of the covering $\overline{S_M}\rightarrow S_M$ is isomorphic to the cyclic group $G_M/A_M$. Denote by $g$ its generator induced by $g_0\in G_M$. Denote by $f$ the birational transformation $\Hol(g_0)$ of $Y$. We have $f\circ \overline{\Dev}=\overline{\Dev}\circ g$. By blowing up $Y$ at some of the indeterminacy points of $f$ (and their infinitely near points), we can and will assume that $f$ is algebraically stable (see \cite{DF01}). The only effect of doing so is to add some extra points into the indeteminacy set $I$ of $\overline{\Dev}$. 

\begin{lem}\label{contractedcurvesareaway}
The contracted curves of $f^n,n\in\ZZ$ are disjoint from the image of $\overline{\Dev}$. 
\end{lem} 
\begin{pro}
Suppose by contradiction that a curve $C\subset Y$ contracted by $f^n$ intersects the image of $\overline{\Dev}$. Since $\overline{\Dev}$ is locally biholomorphic where it is defined, the inverse image $\overline{\Dev}^{-1}(C)$ is a curve on $\overline{S_M}$. Using the relation $f^n\circ \overline{\Dev}=\overline{\Dev}\circ g^n$, we see that $g^n(\overline{\Dev}^{-1}(C))$ is a curve on $\overline{S_M}$ contracted by $\overline{\Dev}$. This is a contradiction as we are already in the case where there are no contracted curves. 
\end{pro}

\begin{lem}
The birational transformation $f=\Hol(g_0)$ is loxodromic.
\end{lem}
\begin{pro}
Suppose by contradiction that $f$ is not loxodromic. We first claim that $f$ preserves a pencil of curves. By definition a \Jonqui twist or a Halphen twist preserves a pencil of curves. Thus we assume that $f$ is elliptic. An elliptic element comes from a holomorphic vector field on $Y$. An elliptic element of infinite order exists only if $Y$ is a rational surface, a ruled surface, an elliptic surface or birational to a surface of Kodaira dimension zero covered by an abelian surface. 

If $Y$ is birational to a surface covered by an abelian surface, then $f$ preserves a transversely Euclidean foliation coming from a linear foliation on the abelian surface. This foliation can be pulled-back by $\overline{\Dev}$ to a foliation on $\overline{S_M}$ invariant under $g$. This further induces a holomorphic foliation on $S_M$ which by Lemma \ref{theonlyfoliations} coincides with one of the two obvious foliations on $S_M$. However neither of these foliations is transversely Euclidean by Lemma \ref{nomeasures}, contradiction.

If $Y$ is rational, then $f$ preserves a pencil of rational curves by Proposition 2.3 of \cite{BD15}. If $Y$ is an elliptic surface of Kodaira dimension one, then $f$ preserves the elliptic fibration of $Y$. If $Y$ is a non-rational ruled surface, then $f$ preserves the rational fibration. The claim follows.

Now we know that $f$ preserves a pencil of curves. This pencil gives rise to a possibly singular holomorphic foliation $S_M$ which by Lemma \ref{theonlyfoliations} coincides with one of the two obvious foliations on $S_M$. By abuse of notation, we use the same letter $\FF$ to denote this foliation on $S_M$ and the one on $\overline{S_M}$. The fact that $\FF$ is induced by a pencil of curves on $Y$ implies that the images of the leaves of $\FF$ by $\overline{\Dev}$ are contained in algebraic curves. The actions of $A_M$ on the spaces of leaves of both of the two foliations are non-discrete, thus the leaves of the two foliations in $\overline{S_M}$ are not closed. Hence the image of a leaf of $\FF$ by $\overline{\Dev}$ cannot be contained in an algebraic curve.
\end{pro}

In the sequel we fix a \Kah{} metric on $Y$ and we endow $\overline{S_M}\backslash I$ with the \Kah{} metric pulled back from $Y$ by $\overline{\Dev}$. Before we consider Ahlfors-Nevanlinna currents, a few words need to be said about the \Kah{} metric. In Section \ref{AhlforsNevanlinnasection}, for constructing Ahlfors-Nevanlinna currents the \Kah{} surface needs to be compact so that the difference between any Riemannian metric and the \Kah{} metric is everywhere bounded by a constant. In our situation here, though $\overline{S_M}\backslash I$ is not compact, we will be able to use freely all the results of Section \ref{AhlforsNevanlinnasection} because of the following three observations: 1) the \Kah{} metric on $\overline{S_M}\backslash I$ is pulled back from the compact surface $Y$; 2) in a small neighborhood $U$ of a point $e\in I$, the map $\overline{\Dev}\vert_{U\backslash\{e\}}$ factorizes through a compact surface (Zariski's factorization theorem for birational maps); 3) the entire curves with which we will deal lie in a compact subset of $\overline{S_M}$. Roughly speaking, these three observations allow us to think of the part of $\overline{S_M}$ on which we will work as an open set of a compact \Kah{} surface.

The set of points in $\overline{S_M}$ that are mapped by $\overline{\Dev}$ to indeterminacy points of $f$ is discrete and countable. The indeterminacy set $I$ of $\overline{\Dev}$ is also discrete and countable. Therefore we can find two fibers $F_a,F_b$ of $\rho:\overline{S_M}\rightarrow \RR^{*+}$ such that:
\begin{itemize}
	\item $F_a\cup F_b$ is disjoint from $I$ and $\overline{\Dev}(F_a\cup F_b)$ is disjoint from the indeteminacy points of $f$;
	\item $g(F_a)=F_b$.
\end{itemize}
We will view the covering map $\HH\times\CC\rightarrow \overline{S_M}$ and the developing map $\Dev:\HH\times\CC\dashrightarrow Y$ (where it is defined) as families of entire curves. By choosing an appropriate path in $\HH$ from a point of vertical coordinate $a$ to a point of vertical coordinate $b$, we can extract from the above family a family of entire curves $(\xi_t)_{t\in [a,b]}$ on $\overline{S_M}$ parametrized by the interval $[a,b]$ such that 
\begin{itemize}
	\item $\forall t\in [a,b]$, the image of $\xi_t:\CC\rightarrow \overline{S_M}$ is disjoint from the indeterminacy set $I$ of $\overline{\Dev}$;
	\item $\xi_t$ parametrizes a leaf of $F_t$; in particular $\xi_a$ (resp. $\xi_b$) parametrizes a leaf of $F_a$ (resp. $F_b$).
\end{itemize} 
We can push the family $(\xi_t)_t$ forward by $\overline{\Dev}$ to obtain a family of entire curves $(\overline{\Dev}\circ \xi_t)_t$ on $Y$. As the covering map and the developing map (where it is defined) are locally biholomorphic, the derivative $\xi_t'(z)$ is non-zero for all $t\in[a,b]$ and for all $z\in\CC$. We claim that the families $(\xi_t)_t$ and $(\overline{\Dev}\circ \xi_t)_t$ are uniform in the sense of Definition \ref{defuniformfamily}. This is clear if $\overline{\Dev}$ has no indeterminacy points because in that case the entire curves $\overline{\Dev}\circ \xi_t$ factorize through the compact sets $F_t$. Since $\overline{\Dev}$ is locally birational, the same reasonning works after blowing up the indeterminacy points contained in the $F_t,t\in [a,b]$. 

By Lemma \ref{lemfamilyofentirecurves}, we can construct a family of Ahlfors-Nevanlinna currents $(\Phi_t)_t$ associated with the uniform family of entire curves $(\xi_t)_t$, after fixing an appropriate sequence of radii once and for all. We construct corresponding Ahlfors-Nevanlinna currents associated with the $\overline{\Dev}\circ \xi_t$: they are the push-forward $\overline{\Dev}_*\Phi_t$. Lemma \ref{lemfamilyofentirecurves} tells us that the cohomology classes $[\overline{\Dev}_*\Phi_t]\in \operatorname{H}^{1,1}(Y,\RR)$ are all the same. We also know that they are nef (see Lemma \ref{ANcurrentisnef}). 

As (the images of) the entire curves $\overline{\Dev}\circ \xi_a$ and $\overline{\Dev}\circ \xi_b$ are disjoint from the contracted curves and the indeterminacy set of $f$ by Lemma \ref{contractedcurvesareaway}, we can push forward the Ahlfors-Nevanlinna current $\overline{\Dev}_*\Phi_a$ by $f$ without any ambiguity. We want to compare the pushed forward current $f_*(\overline{\Dev}_*\Phi_a)$ with $\overline{\Dev}_*\Phi_b$. We have 
\[f_*(\overline{\Dev}_*\Phi_a)=\overline{\Dev}_*(g_*\Phi_a).\]
Thus we just need to compare $g_*\Phi_a$ and $\Phi_b$. By Plante's Theorem \ref{Plantethm}, the closed positive currents $\Phi_a,\Phi_b$ are respectively directed by the laminations $F_a,F_b$. As $g$ sends $F_a$ to $F_b$ preserving their lamination structures, the push forward $g_*\Phi_a$ is a closed positive current directed by $F_b$. By Sullivan's Theorem \ref{Sullivanthm}, the two currents $g_*\Phi_a$ and $\Phi_b$ correspond to two transverse invariant measures on $F_b$. However by Lemma \ref{theonlymeasure}, there exists only one transverse invariant measure on $F_b$ up to multiples. Thus, we have
\[\lambda g_*\Phi_a=\Phi_b \quad \text{for some} \ \lambda\in\RR^{*+}.\]
It follows that
\[\lambda f_*(\overline{\Dev}_*\Phi_a)=\overline{\Dev}_*\Phi_b.\]
By the equality of cohomology classes $[\overline{\Dev}_*\Phi_a]=[\overline{\Dev}_*\Phi_b]$, we get
\begin{equation}\lambda f_*[\overline{\Dev}_*\Phi_a]=[\overline{\Dev}_*\Phi_a].\label{eq:loxodromiceigenclass}\end{equation}
As (the images of) the entire curves $\overline{\Dev}\circ \xi_a$ and $\overline{\Dev}\circ \xi_b$ are disjoint from the contracted curves and the indeterminacy set of $f$ by Lemma \ref{contractedcurvesareaway}, we get also
\begin{equation}f^*[\overline{\Dev}_*\Phi_a]=\lambda[\overline{\Dev}_*\Phi_a].\label{eq:loxodromiceigenclass2}\end{equation}

\begin{lem}
The dynamical degree of $f$ is equal to $\lambda$ or $\lambda^{-1}$. 
\end{lem}
\begin{pro}
Denote the dynamical degree of $f$ by $\lambda(f)$. There exists a unique nef cohomology class $v_f^+$ such that $f^*v_f^+=\lambda(f)v_f^+$. By Proposition 1.11 of \cite{DF01}, we have the following equality for intersection numbers:
\begin{equation}
(f^*v_f^+,[\overline{\Dev}_*\Phi_a])=(v_f^+,f_*[\overline{\Dev}_*\Phi_a]).\label{eq:pushpullformula}
\end{equation}
If $[\overline{\Dev}_*\Phi_a]$ and $v_f^+$ are proportional, then $\lambda(f)=\lambda$ by Equation \eqref{eq:loxodromiceigenclass2}. Assume that they are not proportional; this implies that their intersection is strictly positive because they are both nef. Then Equations \eqref{eq:loxodromiceigenclass} and \eqref{eq:pushpullformula} force the equality between $\lambda(f)$ and $1/\lambda$.
\end{pro}

Replacing $f$ with $f^{-1}$ if necessary, we can and will assume that the dynamical degree of $f$ is $\lambda$.

\begin{lem}
We can assume that $f$ acts by automorphism on $Y$, without loosing any other property that we need.
\end{lem}
\begin{pro}
We first prove that all the irreducible curves contracted by the iterates $f^n$ are of strictly negative self-intersection. Let $E$ be an irreducible curve contracted by $f^n$. By Lemma \ref{contractedcurvesareaway}, the curve $E$ is disjoint from $\overline{\Dev}(F_a)$ which is the support of $\overline{\Dev}_*\Phi_a$. Therefore the intersection number $[\overline{\Dev}_*\Phi_a]\cdot E$ is zero. As $[\overline{\Dev}_*\Phi_a]$ is nef, we have $[\overline{\Dev}_*\Phi_a]^2\geq 0$. It follows from Hodge index theorem that $E^2\leq 0$, with equality if and only if $[\overline{\Dev}_*\Phi_a]^2=0$ and $E$ is proportional to $[\overline{\Dev}_*\Phi_a]$. Since $[\overline{\Dev}_*\Phi_a]$ is an eigenvector associated with $\lambda$, the equality $[\overline{\Dev}_*\Phi_a]^2=0$ would imply that the algebraically stable map $f$ is an automorphism (see Section \ref{birpre}). But then $[\overline{\Dev}_*\Phi_a]$ would be irrational and could not be proportional to $E$.  

We write the Zariski factorization of $f$ as $Y\leftarrow \hat{Y}\rightarrow Y$. Let $E_1,\cdots,E_m$ be the irreducible curves contracted by $f$. Denote by $\hat{E}_1,\cdots,\hat{E}_m$ their strict transforms in $\hat{Y}$. Among the $\hat{E}_i$, there exists at least one $(-1)$-curve, let us say, $\hat{E}_1$. Since $\hat{Y}$ is obtained from $Y$ by blow-ups, we have $\hat{E}_1^2\leq E_1^2$. We have showed that $E_1^2<0$. It follows that $E_1$ is already a $(-1)$-curve on $Y$. Now we contract it to obtain a new surface $Y_1$. We need to verify that all the hypothesis still hold on $Y_1$. The contraction may give rise to new curves contracted by $f$, but the new contracted curves on $Y_1$ come from the curves on $Y$ contracted by $f^2$. So they are still disjoint from the image of $\overline{\Dev}$ and are of strictly negative self-intersection on $Y_1$. Hence we can continue the process. This process terminates because the Picard number drops down by one after each step. At last we get a surface on which $f$ contracts no curves, i.e.\ $f$ acts by automorphism. 
\end{pro}

Once we know that $f$ is a loxodromic automorphism, Theorem \ref{CantatDinhSibony} implies that $\overline{\Dev}_*\Phi_a$ is the unique closed positive current with cohomology class $[\overline{\Dev}_*\Phi_a]$. However the cohomology class of $\overline{\Dev}_*\Phi_b$ is also $[\overline{\Dev}_*\Phi_a]$. This leads to a contradiction because $\overline{\Dev}_*\Phi_a$ and $\overline{\Dev}_*\Phi_b$ are two different currents. Indeed their supports are respectively $\overline{\Dev}(F_a), \overline{\Dev}(F_b)$ and $\overline{\Dev}(F_a)\neq\overline{\Dev}(F_b)$ because otherwise $\overline{\Dev}$ would induce a map from $S_M$ to $Y$. The proof of Proposition \ref{propnoncyclic} is finished.

\begin{rmk}
The very existence of immersed Levi-flat hypersurfaces as $\overline{\Dev}(F_a)$ imposes strong restrictions on the geometry of $Y$. For example there are no such immersed Levi-flat hypersurfaces in $\PP^2$ (see \cite{Der05}). However there exist families of Levi-flat hypersurfaces on other surfaces and we are not able to conclude directly by the existence of an "immersion" of $\overline{S_M}$ into $Y$. This is why the geometry of the cyclic covering $\overline{S_M}\rightarrow S_M$ plays a crucial role in our proof.
\end{rmk}

\subsubsection{Injective holonomy}\label{subsectioninjectiveholonomy}

Since we have proved that the image of the holonomy is not cyclic, its kernel must be finite by Lemma \ref{kernellem}. Thus changing $S_M$ into a finite cover \emph{we can and will assume in the sequel that $\Hol$ is injective.} We will identify $G_M$ with its image $\Hol(G_M)\in\Bir(Y)$. 

\begin{lem}\label{keylem1}
The group $G_M$ is an elliptic subgroup of $\Bir(Y)$.
\end{lem}
\begin{pro}
We apply Theorem \ref{Deserti} to the solvable group $G_M\subset \Bir(Y)$. Up to conjugating the holonomy representation there are five possibilities in Theorem \ref{Deserti} and we need to rule out the last four ones.

In case 5) $Y$ is an abelian surface, the group $G_M$ is generated by translations and a loxodromic automorphism. The stable and the unstable foliations of the loxodromic automorphism (see Example 1.1 of \cite{CF03}), which are both linear foliations on $Y$, are preserved by $G_M$. Thus they can be pulled back to two holomorphic foliations on $S_M$. Induced by linear foliations on $Y$, these two pulled back foliations are transversely Euclidean. But they must coincide with the obvious foliations on $S_M$ by Lemma \ref{theonlyfoliations}; this contradicts Lemma \ref{nomeasures}.

In case 4) $Y$ is rational and $G_M$ is in $\CRC$. In this case $A_M$ is contained in $\{(\alpha x,\beta y)\vert \alpha,\beta\in \CC^*\}$ and $g_0$ is a monomial map $(x^py^q,x^ry^s)$ such that the matrix $B=\begin{pmatrix}p&q\\r&s\end{pmatrix}\in \GL_2(\ZZ)$ is hyperbolic, i.e.\ has one eigenvalue $>1$ and one $<1$. The conjugation action of $g_0$ on $A_M$ is given by $(\alpha,\beta)\mapsto (\alpha^p\beta^q,\alpha^r\beta^s)$. The exponential map semi-conjugates the action of $B=\begin{pmatrix}p&q\\r&s\end{pmatrix}$ on $\CC^2$ to the action of $g_0$ on $\CC^*\times \CC^*$. We think of $A_M$ as a $\ZZ$-module of rank $3$ with an irreducible action of $g_0$. Its preimage $\overline{A_M}$ in $\CC^2$ by the exponential map is a $\ZZ$-module of rank $5$ invariant under $B$. The kernel of the exponential map, generated by $(0,2\pi i)$ and $(2\pi i,0)$ is invariant under $B$ and the action is irreducible. Hence either $\overline{A_M}$ is an indecomposable module of rank $5$ or there is an indecomposable submodule of rank $3$ which is isomorphic to $A_M$. However an indecomposable module, subgroup of $\CC^2$, has to be of even rank because $B$ is a hyperbolic matrix. We obtain thus a contradiction. Hence case 4) is not possible.

Case 3) of Theorem \ref{Deserti} is impossible because $G_M$ is not virtually abelian.

It suffices to show that case 2) of Theorem \ref{Deserti} is not possible for $G_M$. Suppose the contrary. The rational fibration preserved by $G_M$ can be pulled-back to a holomorphic foliation on $\HH\times\CC$. By the equivariance of $D$, this equips $S_M$ with a holomorphic foliation. This foliation must coincide with one of the two obvious foliations on $S_M$ by Lemma \ref{theonlyfoliations}. Acting on $\HH\times\CC$, the elements $g_0,\cdots,g_3$ permute the leaves of the foliation. The action of $G_M$ on the spaces of leaves of the two foliations on $\HH\times\CC$, i.e.\ its actions on the $\CC$-factor and on the $\HH$-factor are both non-discrete. This means that, on the $\Bir(Y)$ side, the action of $G_M$ on the base of the rational fibration is non discrete. As the automorphism group of a curve of general type is finite, the base is either $\PP^1$ or an elliptic curve. But the base cannot be an elliptic curve neither because otherwise the fibration would be transversely Euclidean. Thus the base of the rational fibration is $\PP^1$ and $Y$ is a rational surface. 

We have a morphism $\sigma:G_M\rightarrow \PGL_2(\CC)$ that records the action of $G_M$ on the base of the rational fibration. Since this base action is non-discrete, Lemma \ref{kernellem} implies that $\sigma$ is injective. So $g_1,g_2,g_3$ have infinite actions on the base. Suppose by contradiction that one of the $g_i$, say $g_1$, is a \Jonqui twist. Theorem \ref{centralizer} case 3 says that for any $h$ that commutes with $g_1$, the actions of $h$ and $g_1$ on the base generate a virtually cyclic group. But we said that the actions of $g_1,g_2$ and $g_3$ on the base generate a group isomorphic to $\ZZ^3$, contradiction. Hence $g_1,g_2,g_3$ are all elliptic elements of $\Bir(Y)$ and $A_M$ is an elliptic subgroup of $\Bir(Y)$. Up to replacing $A_M$ by a finite index free abelian subgroup, we can assume that $A_M$ is contained in $\Aut^0(Z)$, the connected component of the automorphism group of a rational surface $Z$. The group $\Aut^0(Z)$ is an algebraic group; we denote by $\overline{A_M}$ the Zariski closure of $A_M$ in $\Aut^0(Z)$. Since $A_M$ is infinite, $\overline{A_M}$ is an algebraic group of dimension $\geq 1$.

We want to prove that no element of $G_M$ is a \Jonqui twist. For this purpose we apply an argument used by S. Cantat in the appendix of \cite{DP12}. Any element of $G_M$ normalizes $A_M$, thus normalizes $\overline{A_M}$. We have two possibilities for the action of the abelian algebraic group $\overline{A_M}$ on the rational surface $Z$, either it has a Zariski open orbit, or its orbits form a pencil of curves. 

Assume that the orbits of $\overline{A_M}$ form a pencil of curves. This pencil of curves must differ from the original rational fibration preserved by $A_M$ because the actions of $g_1,g_2,g_3$ on the base of the rational fibration are infinite. Every element of $G_M$ normalizes $\overline{A_M}$, it preserves this pencil of curves. Recall that every element of $G_M$ preserves also the rational fibration, thus preserves simultaneously two pencils of curves. This implies that $g_0$ is an elliptic element (cf. \cite{DF01}). Therefore the group $G_M$ contains no \Jonqui twists.

Now assume that $\overline{A_M}$ has a Zariski open orbit $O$. We have three possibilities for $O$; it is a principal homogeneous space isomorphic to $\CC^2$, $\CC\times\CC^*$ or $\CC^*\times\CC^*$. Since an element of $G_M$ normalizes $\overline{A_M}$, it acts on $O$ by automorphism of principal homogeneous spaces. If $O=\CC^2$ then every element of $G_M$ would be affine, thus elliptic. If $O=\CC\times\CC^*$ then an element of $G_M$ would be of the form $(a x +b,\alpha y)$ with $a,\alpha\in\CC^*,b\in\CC$, which is again elliptic. If $O=\CC^*\times\CC^*$ then every element of $G_M$ would be contained in the group generated by $\{(\alpha x,\beta y)\vert \alpha,\beta \in\CC^*\}$ and $\{(x^py^q,x^ry^s),\begin{pmatrix}p&q\\r&s\end{pmatrix}\in \operatorname{GL}_2(\ZZ)\}$. In this case an element of $G_M$ is either elliptic or loxodromic, but it cannot be loxodromic because we work already under the hypothesis that $G_M$ preserves a rational fibration. Thus we have proved that every element of $G_M$ is elliptic. This implies that $G_M$ is an elliptic subgroup by Theorem \ref{Deserti}. 
\end{pro}

We proved that $G_M$ is an elliptic subgroup of $\Bir(Y)$. Up to taking a finite index subgroup, $G_M$ is contained in $\Aut^0(Z)$, the component of identity of the automorphism group of a projective surface $Z$ birational to $Y$. The $\Aut^0$ of a projective variety is an algebraic group. By Chevalley structure theorem it is an extension of an abelian variety by a linear algebraic group. As $\Aut^0(Z)$ contains the non-abelian infinite group $G_M$, its linear part is not trivial. Hence $Z$ is ruled by Theorem 14.1 of \cite{Ueno75} (see also \cite{Bru00} Chapter 6.3). If $Z$ were a non-rational ruled surface, then $G_M$ would preserve the ruling and the ruling would be pulled back by $\Dev$ to one of the two obvious foliations on $S_M$. Using the fact that $G_M$ acts non-discretely on the space of leaves and the fact that the two obvious foliations are not transversely Euclidean, we obtain a contradiction as in the proof of Lemma \ref{keylem1}. 

Therefore $Z$ is a rational surface. Since $G_M$ is solvable, it comes from a group of automorphisms of a Hirzebruch surface. \emph{We can and will assume in the sequel that $Z=\FFh_n$ is a Hirzebruch surface and that $G_M\subset \Aut^0(Z)$} (cf. Corollary \ref{switchcor1}). Note that from now on we take $Z$ as the target space of the developing map.

Since $\Aut(\FFh_n)$ preserves the rational fibration on $\FFh_n$, we have a group homomorphism $\sigma:G_M\rightarrow \PGL_2(\CC)$ which encodes the action of $G_M$ on the base $\PP^1$ of the rational fibration. As $G_M$ is solvable, we can assume, maybe after replacing $G_M$ with a subgroup of index two, that $\sigma(G_M)\subset \PGL_2(\CC)$ fixes at least one point in $\PP^1$. Let us decompose $\PP^1$ as $\CC\cup\{\infty\}$ where $\infty$ is one of the fixed point of $\sigma(G_M)$. As in the proof of Lemma \ref{keylem1}, the rational fibration induces a foliation on $S_M$ which must coincide with one of the two obvious foliations on $S_M$; and we deduce from this that $\sigma(A_M)$ is not discrete. By Lemma \ref{kernellem}, this implies that $\sigma$ is injective (up to taking a finite index subgroup). By Lemma \ref{injectionintoPGLtC}, we can write $\sigma(g_1),\sigma(g_2),\sigma(g_3)$ as $x\mapsto x+u_i$ for some $u_i\neq 0$, and $\sigma(g_0)$ as $x\mapsto \nu x$ for some $\nu\in\CC^*$ of infinite order.

\begin{lem}\label{developingmapbiholomorphic}
The developing map $\Dev:\HH\times\CC\dashrightarrow \FFh_n$ is everywhere defined and is locally biholomorphic.
\end{lem}
\begin{pro}
First we claim that $\Dev$ contracts no curves. Suppose by contradiction that $\Dev$ contracts a curve $C\subset\HH\times\CC$. Let $\gamma\in G_M$ be a non-trivial element. From the relation $\Dev\circ\gamma=\gamma\circ\Dev$ and the fact that $\gamma$ acts by automorphism on $\FFh_n$, we deduce that $\gamma(C)\subset \HH\times\CC$ is also a contracted curve of $\Dev$. Since locally there is only a finite number of contracted curves, the union $\bigcup_{\gamma\in G_M}\gamma(C)$ is a $G_M$-invariant set locally closed in $\HH\times\CC$. Therefore the image of $C$ in the quotient $S_M$ is locally closed, i.e.\ it is a curve on $S_M$. This contradicts the fact that $S_M$ has no curves.

Suppose by contradiction that $p\in\HH\times\CC$ is an indeterminacy point of $\Dev$. Take a local chart of birational structure $U$ at $p$. We factorize $\Dev\vert_U:U\dashrightarrow \FFh_n$ as $U\xleftarrow{\pi_1} V\xrightarrow{\pi_2} \FFh_n$ where $\pi_1$ is a composition of (inverses of) blow-ups at $p$ and its infinitely near points, and $\pi_2$ is an open embedding because $\Dev$ contracts no curves. Note that here we have holomorphic foliations on $U$ and $V$, pulled back from the rational fibration on $\FFh_n$. As the foliation on $U$ is regular, the exceptional curve of $\pi_1$ is an invariant curve of the foliation on $V$. This implies that its image by $\pi_2$ into $\FFh_n$ must be contained in a fiber of the rational fibration. However on $\FFh_n$ there are no $(-1)$-curves contained in the fibers of the rational fibration, this contradicts the fact that $\pi_2$ is an open embedding. 
\end{pro}

The above lemma tells us that the birational structure on $S_M$ is in fact a $(\Aut(\FFh_n),\FFh_n)$-structure in the classical sense. 

\begin{lem}\label{noinvcurves}
Any $G_M$-invariant curve is disjoint from the image of $\Dev$.
\end{lem}
\begin{pro}
Let $C$ be a curve which intersects the image of $\Dev$, then $\Dev^{-1}(C)$ is a curve on $\HH\times\CC$. Note that elements of $G_M$ are regular on the intersection of $C$ with the image of $\Dev$. So if $C$ were $G_M$-invariant, then $\pi(\Dev^{-1}(C))$ would be a curve on $S_M$.
\end{pro}

\textit{First case: Assume that the Hirzebruch surface $\FFh_n$ is not $\PP^1\times\PP^1$, i.e.\ $n\geq 1$.} The fiber over $\infty\in\PP^1$ and the exceptional section of $\FFh_n$ are $G_M$-invariant curves. Lemma \ref{noinvcurves} implies that the image of the developing map is contained in the complement of these two invariant curves which is isomorphic to $\CC^2$. The automorphisms $g_1,g_2,g_3$ are of the form $(x,y)\mapsto (x+u_i,\beta_i y+R_i(x))$ where $\beta_i\in\CC^*$ and $R_i$ is a polynomial of degree $\leq n$; the automorphism $g_0$ is $(x,y)\mapsto (\nu x,\beta_0 y+R_0(x))$ where $\beta_0\in\CC^*$ and $R_0$ is a polynomial of degree $\leq n$. 

Assume first that $\beta_1,\beta_2,\beta_3$ are not all equal to $1$, for example $\beta_1\neq 1$. In this case $g_1$ has a fixed point $e$ on the fiber over $\infty$ which is not on the exceptional section of $\FFh_n$. By commutativity, the point $e$ is fixed by $A_M$; then by the fact that $g_0$ normalizes $A_M$, the whole group $G_M$ fixes $e$. We can blow-up $e$ and contract the strict transform of the initial fiber to get $\FFh_{n-1}$. The group $G_M$ remains a group of automorphisms of $\FFh_{n-1}$. Moreover the image of the developing map is not affected by this elementary transformation. Therefore the initial birational structure reduces to a $(\Aut(\FFh_{n-1}),\FFh_{n-1})$-structure. We continue this process and reduce the birational structure to a $(\Aut(\FFh_{1}),\FFh_{1})$-structure. The Hirzebruch surface $\FFh_1$ is the blow-up of $\PP^2$ at one point; the exceptional divisor is the exceptional section and is disjoint from the image of $\Dev$. Therefore we finally get a $(\PGL_3(\CC),\PP^2)$-structure.

Assume that $\beta_1=\beta_2=\beta_3=1$. We can conjugate $g_1$ inside $\CRC$, by elements of the form $(x,y)\dashrightarrow (x,y+\delta x^d)$, to decrease the degree of $R_1$ until $g_1$ becomes $(x,y)\mapsto (x+u_1,y)$; note that this only modifies the fiber at $\infty$ so that the conjugation does not affect $\Dev$. After these conjugations, $g_2,g_3$ become $(x,y)\dashrightarrow (x+u_i,y+\tilde{R}_i(x)),\ i=2,3$ and $g_0$ becomes $(x,y)\dashrightarrow (\nu x,\beta_0 y+\tilde{R}_0(x))$ where the $\tilde{R}_i$ are polynomials for $i=0,2,3$. The commutation relations between $g_1$ and $g_2,g_3$ reads:
\[
\tilde{R}_i(x)=\tilde{R}_i(x+u_1), \quad i=2,3;
\]
this implies immediately that $\tilde{R}_2$ and $\tilde{R}_3$ are constants. Therefore we have conjugated $A_M$ to a subgroup of $\PGL _3(\CC)$. Now the transformation $g_0\circ g_1\circ g_0^{-1}$ is 
\[(x,y)\dashrightarrow (x+\nu u_1,y+\tilde{R}_0(\nu^{-1}x+u_1)-\tilde{R}_0(\nu^{-1}x))\]
For $g_0\circ g_1\circ g_0^{-1}$ to be in $A_M$, the polynomial function $\tilde{R}_0(\nu^{-1}x+u_1)-\tilde{R}_0(\nu^{-1}x)$ needs to be a constant. This implies that the degree of $\tilde{R}_0$ is at most $1$, i.e.\ $g_0$ is also in $\PGL _3(\CC)$. We get again a $(\PGL_3(\CC),\PP^2)$-structure.

\textit{Second case: the Hirzebruch surface is $\PP^1\times \PP^1$.} Considering a finite unramified cover of $S_M$, we can assume that $G_M$ is included in the identity component of the automorphism group which is $\PGL _2(\CC)\times \PGL _2(\CC)$. Replacing $G_M$ with a index two subgroup if necessary, we have two injective homomorphisms $\sigma_1,\sigma_2$ from the solvable group $G_M$ to $\PGL _2(\CC)$. The image $\sigma_1(G_M)$ (resp. $\sigma_2(G_M)$) fixes at least one point in the first (resp. second) factor $\PP^1$. Removing the two corresponding $G_M$-invariant curves from $\PP^1\times\PP^1$, we get a Zariski open set which is isomorphic to $\CC^2$ and in which the image of the developing map is contained. This means that the birational structure is reduced to a complex affine structure. 

Bruno Klingler proved in \cite{Kli98} that the only $(\PGL_3(\CC),\PP^2)$-structure on $S_M$ is the natural one, this finishes the proof of Theorem \ref{mainthm} for Inoue surfaces of type $S^0$.

\section{Inoue surfaces of type $S^{\pm}$}\label{thirdsection}
\subsection{Description}
Let $n\in\NN^*$. Consider the group of upper-triangular matrices
\[
\Lambda_n=\left\{ \begin{pmatrix}1&x&\frac{z}{n}\\0&1&y\\0&0&1\end{pmatrix}, \quad x,y,z\in\ZZ \right\}.
\]
The center of $\Lambda_n$ is the infinite cyclic group $C_n$ generated by $\begin{pmatrix}1&0&\frac{1}{n}\\0&1&0\\0&0&1\end{pmatrix}$. The quotient $\Lambda_n/C_n$ is isomorphic to $\ZZ^2$. Let $N\in \operatorname{SL} _2(\ZZ)$ be a matrix with eigenvalues $\alpha,\frac{1}{\alpha}$ such that $\alpha>1$. Let $\varphi$ be an automorphism of the group of real upper-triangular matrices which preserves $\Lambda_n$, acts trivially on $C_n$ and acts on $\Lambda_n/C_n\cong \ZZ^2$ as $N$. We form a semi-direct product $\Gamma_N=\Lambda_n\rtimes \ZZ$ where the $\ZZ$ factor acts on $\Lambda_n$ as $\varphi$. The group $\Gamma_N$ acts on the group of real upper-triangular matrices which is identified with $\RR^3=\RR\times\CC$. Define an action of $\Gamma_N$ on $\HH\times\CC=\RR^{>0}\times\RR\times\CC$ with $\Lambda_n$ acting trivially on $\RR^{>0}$ and $1\in\ZZ$ acting on $\HH$ as $x\mapsto \alpha x$. This action is holomorphic and the quotient $S_N=\HH\times\CC/\Gamma_N$ is a compact non-\Kah{} surface called an \emph{Inoue surface of type $S^+$} (\cite{Ino74}). Note that the Inoue surface depends on $n,\varphi$, and $\varphi$ depends on $N$; we denote it by $S_N$ because $N$ is the most significant parameter. 

The group $\Gamma_N$ can be identified with a lattice in one of the two following solvable Lie groups which are subgroups of $\Aff _2(\CC)$ (cf. \cite{Kli98}):
\[
\operatorname{Sol}^1=\left\{\begin{pmatrix}1 &a&b\\0&d &c\\0&0&1\end{pmatrix},a,b,c,d\in\RR,d>0\right\},\ 
\operatorname{Sol}^{1'}=\left\{\begin{pmatrix}1 &a&b+i\log(d)\\0&d &c\\0&0&1\end{pmatrix},a,b,c,d\in\RR,d>0\right\}.
\]
Conversely any torsion free lattice of these two groups gives an Inoue surface of type $S^+$. Note that a finite unramified cover of an Inoue surface of type $S^+$ is an Inoue surface of type $S^+$. 

Concretely $\Gamma_N$ has four generators $g_0,g_1,g_2,g_3$ which act on $\HH\times \CC$ as:
\begin{align*}
g_0:(x,y)&\mapsto (\alpha x,y+t)\\
g_i:(x,y)&\mapsto (x+a_i,y+b_i x+c_i)\quad i=1,2\\
g_3:(x,y)&\mapsto (x,y+\frac{b_1a_2-b_2a_1}{n})
\end{align*}
where $t$ is a complex number, $(a_1,a_2)$ (resp. $(b_1,b_2)$) is a real eigenvector of $N$ corresponding to the eigenvalue $\alpha$ (resp. $\alpha^{-1}$) and $c_1,c_2$ are some complex numbers (see \cite{Ino74} for the explicit expressions of $c_1,c_2$). The center $C_n$ of $\Lambda_n$ is also the center of $\Gamma_N$, it is generated by $g_3$. The normal subgroup $\Lambda_n$ is generated by $g_1,g_2,g_3$. We have 
\begin{align*}
&g_1^{-1}g_2^{-1}g_1g_2=g_3^n\\
&g_0g_ig_0^{-1}=g_1^{n_{i1}}g_2^{n_{i2}}g_3^{m_i}, \quad i=1,2
\end{align*}
where $n_{ij}$ are entries of the matrix $N$ and $m_1,m_2$ are two integers depending on $c_1,c_2$.

The Inoue surface $S_N$ has an obvious $(\Aff _2(\CC),\CC^2)$-structure. The surface $\overline{S_N}=\HH\times\CC/\Lambda_n$ is an infinite cyclic covering space of $S_N$. As a real manifold, $\overline{S_N}$ admits a fibration $\rho:\overline{S_N}\rightarrow \RR^{*+}$ where the $\RR_+^{*}=\{t\sqrt{-1},t\in\RR_+^*\}$ is the vertical axis of $\HH\subset \HH\times\CC$. The fibers, denoted by $E_t$, are quotients of $\{x+t\sqrt{-1},x\in\RR^{*+}\}\times \CC$ by $\Lambda_n$; they are compact real nilmanifolds of dimension $3$. The $E_t$ are Levi-flat hypersurfaces in $\overline{S_N}$ and they are foliated by entire curves coming from the vertical complex lines in $\HH\times\CC$.

The analogues of Lemmata \ref{kernellem}, \ref{injectionintoPGLtC}, \ref{theonlyfoliations}, \ref{theonlymeasure} and \ref{nomeasures} still hold (we omit the details when the proof is exactly the same).
\begin{lem}\label{kernellem2}
If $K$ is a non-trivial normal subgroup of $\Gamma_N$, then $K$ has finite index in $C_n$, $\Lambda_n$ or $\Gamma_N$.
\end{lem}
\begin{pro}
The conjugation action of $g_0$ on $\Lambda_n/C_n$ is just the action of $N\in\operatorname{SL} _2(\ZZ)$ on $\ZZ^2$; it has no eigenvectors in $\ZZ^2\backslash\{0\}$. Thus, if $K$ contains an element of $\Lambda_n$ which is not in $C_n$, then it contains $\Lambda_n$. To conclude, we need only remark that, by the semi-direct product structure, the intersection of a normal subgroup of $\Gamma_N$ with $\Lambda_n$ cannot be trivial.
\end{pro}

\begin{lem}\label{injectionintoPGLtC2}
Let $\sigma:\Gamma_N\rightarrow \PGL_2(\CC)$ be a morphism whose kernel is $C_n$. Then for some affine coordinate $\PP^1=\{x\in\CC\}\cup \{\infty\}$, the images $\sigma(g_i), i=0,\cdots,2$, viewed as homographies of $\PP^1$, may be written as
\begin{align*}
\sigma(g_i):& x\mapsto x+u_i, \quad i=1,2\\ 
\sigma(g_0):& x\mapsto \nu x
\end{align*}
for some $\nu,u_i\in\CC^*$.
\end{lem}

\begin{lem}\label{theonlyfoliations2}
The only (possibly singular) holomorphic foliation on $S_N$ is the obvious one coming from the vertical foliation of $\HH\times \CC$. 
\end{lem}
\begin{pro}
Let $\FF$ be a foliation on $S_N$. As in the proof of Lemma \ref{theonlyfoliations}, we infer that $\FF$ is saturated and non-singular. We denote by $T$ the tangent bundle of $S_M$, by $T^*$ its dual and by $K$ the canonical bundle of $S_M$. We denote by $F_0$ the normal bundle of the obvious foliation (here we use notations of \cite{Ino74}) and by $F$ the normal bundle of $\FF$. The foliation $\FF$ corresponds to a non-zero global section of $T^*\otimes F$. It is proved in \cite{Ino74} that $T^*\otimes F$ has non-zero sections if and only if $F=F_0$. In other words, $\FF$ and the obvious foliation share the same normal bundle. It is also proved in \cite{Ino74} that the space of global sections of $T^*\otimes F_0$ is one dimensional. Thus, $\FF$ must coincide with the obvious foliation.
\end{pro}

\begin{lem}\label{theonlymeasure2}
Up to multiples, there is only one transverse invariant measure on $E_t$.
\end{lem}
\begin{lem}\label{nomeasures2}
The obvious foliation on $S_N$ is not transversely Euclidean.
\end{lem}

Inoue surfaces of type $S^{-}$ are defined similarly: instead of choosing $N$ in $\operatorname{SL} _2(\ZZ)$, we take a matrix in $\operatorname{GL} _2(\ZZ)$ with determinant $(-1)$. Every Inoue surface of type $S^{-}$ has a double unramified cover which is an Inoue surface of type $S^{+}$. Thus, for our purpose it is sufficient to consider only the Inoue surfaces of type $S^{+}$.

\subsection{Proof of Theorem \ref{mainthm} for Inoue surfaces of type $S^+$}
Many details of the proof will be very similar to the case of Inoue surfaces of type $S^0$; we will make them brief.

Equip $S_N$ with a $(\Bir(X),X)$-structure and let $(Y,\Hol,\Dev)$ be a holonomy/developing triple.

\subsubsection{Run again the previous proof}
Lemma \ref{kernellem2} says that there are only four possibilities for the holonomy representation. It is easy to rule out the first possibility: if the holonomy had finite image then the developing map would induce a meromorphic locally birational map from a finite unramified cover of $S_N$ to $Y$, contradicting the fact that $S_N$ has algebraic dimension zero. 

If the kernel $K$ of the holonomy has finite index in $\Lambda_n$, then $K\rtimes \ZZ$ has finite index in $\Gamma_N$; in this case by considering the corresponding finite unramified cover of $S_N$ and the induced birational structure, we can suppose that $K=\Lambda_n$. The image of the holonomy is then cyclic. This is not possible: Lemmata \ref{theonlyfoliations2}, \ref{theonlymeasure2} and \ref{nomeasures2} ensure that the proof of Subsection \ref{subsectionnoncylicholonomy} works exactly in the same way for $S_N$. 

We now rule out the case where the kernel $K$ of $\Hol$ has finite index in $C_n$; we will examine the situation of injective holonomy in the next subsection. After taking a finite unramified cover of $S_N$, we can and will assume that $K=C_n$. Thus, we have an embedding of $\Omega_N=\Gamma_N/C_n\cong \ZZ^2\rtimes \ZZ$ into $\Bir(Y)$. The situation is almost the same as in the case of Inoue surface of type $S^0$; there we had $\ZZ^3\rtimes \ZZ$, here we have $\ZZ^2\rtimes \ZZ$. We can almost copy the proof of Section \ref{subsectioninjectiveholonomy}; we give here a sketch. 

Firstly we prove as in Lemma \ref{keylem1} that $\Omega_N$ is an elliptic subgroup of $\Bir(Y)$. The only difference in the proof is the fourth case of Theorem \ref{Deserti}. In case 4), $\Omega_N$ is contained in the group generated by $\{(\alpha x,\beta y)\vert \alpha,\beta \in\CC^*\}$ and one monomial transformation $(x^py^q,x^ry^s)$ where $\begin{pmatrix}p&q\\r&s\end{pmatrix}\in \operatorname{GL}_2(\ZZ)$. In this case $\Omega_N$ preserves two holomorphic foliations defined by $\iota_1xdy+\nu_1ydx$ and $\iota_2xdy+\nu_2ydx$ where $(\iota_i,\nu_i)\ i=1,2$ are two eigenvectors of $\begin{pmatrix}p&r\\q&s\end{pmatrix}$. These two $\Gamma_N$-invariant foliations induce two foliations on $S_N$; this is impossible because there exists only one holomorphic foliation on an Inoue surface of type $S^+$ by Lemma \ref{theonlyfoliations2}. 

Once we know that $\Omega_N$ is an elliptic subgroup, we prove as in Section \ref{subsectioninjectiveholonomy} that the $(\Bir(X),X)$-structure is reduced to a $(\Aut(\FFh_n),\FFh_n)$-structure, and then to a $(\PGL_3(\CC),\PP^2)$-structure; the arguments here and there are exactly the same. However the only $(\PGL_3(\CC),\PP^2)$-structure on $S_N$ is the obvious one by \cite{Kli98} and its holonomy is injective, a contradiction to the hypothesis that the kernel of the holonomy is $C_n$. Thus, we have proved:
\begin{lem}\label{injectiveholonomy2}
The kernel of the holonomy representation $\Hol$ is trivial.
\end{lem}

\subsubsection{Injective holonomy}
After Lemma \ref{injectiveholonomy2} we know that the holonomy representation is injective. From now on we consider $\Gamma_N$ as a subgroup of $\Bir(Y)$. We apply Theorem \ref{Deserti} to $\Gamma_N$. Case 5) is not possible because the stable and the unstable foliations of a loxodromic automorphism on an abelian surface would induce two transversely Euclidean foliations on $S_N$. Case 4) is not possible because the derived length of $\CC^*\times\CC^*\rtimes \langle\text{a monomial map}\rangle$ is $2$ and that of $\Gamma_N$ is $3$ (here we can also use the foliation argument). Case 3) is impossible because $\Gamma_N$ is not virtually abelian. The following lemma says that case 2) is not possible either.
\begin{lem}
If the group $\Gamma_N$ preserves a rational fibration, then it contains no \Jonqui twists and $Y$ is rational.
\end{lem}
\begin{pro}
The rational fibration preserved by $\Gamma_N$ induces a holomorphic foliation on $S_N$ which coincides with the natural one. The action of $g_3$, even on $\HH\times\CC$, does not permute the leaves, so its action on the base of the rational fibration must be trivial. As regards the action of $\Gamma_N\backslash C_n$ on the base, it is non-discrete by considering the action on the space of leaves. Together with the fact that the foliation is not transversely Euclidean, this implies that the base of the rational fibration is necessarily $\PP^1$. Thus $Y$ is a rational surface.  

Using again the non-discreteness of the base action, we have an embedding $\sigma:\Omega_N=\Gamma_n/C_n\rightarrow\PGL_2(\CC)$. By Lemma \ref{injectionintoPGLtC2}, we infer that $\sigma(g_0),\sigma(g_1),\sigma(g_2)$ are respectively $x\mapsto \gamma x,x\mapsto x+u_1,x\mapsto x+u_2$ where $\gamma,u_1,u_2\in\CC^*$ are such that $\gamma u_i=n_{i1}u_1+n_{i2}u_2$ for $i=1,2$. The sequel of the proof is purely about the group of birational transformations.

Every element of $\Gamma_N$ commutes with $g_3$; for $i=0,1,2$, the group generated by $g_i,g_3$ is isomorphic to $\ZZ^2$. By Theorem \ref{centralizer}, $g_3$ must be an elliptic element. Up to conjugation, $g_3$ is $(x,y)\mapsto (x,y+v_3)$ or $(x,y)\mapsto (x,\nu y)$. Let us first suppose that $g_3$ is $(x,y)\mapsto (x,\nu y)$ for some $\nu\in\CC^*$ of infinite order. By Theorem \ref{ellipticthm} $g_0,g_1,g_2$ are respectively $(\gamma x, R_0(x)y)$ and $(x+u_i,R_i(x)y),\ i=1,2$ where $R_0,R_1,R_2\in\CC(x)$. The relation $g_1^{-1}g_2^{-1}g_1g_2=g_3^n$ reads 
\[
R_2(x)R_1(x+u_2)R_2(x+u_1)^{-1}R_1(x)^{-1}=\nu^n.
\]
For $i=1,2$ write $R_i$ as $\frac{P_i}{Q_i}$ with $P_i,Q_i\in \CC[x]$. Then the above equation becomes
\[
\frac{P_2(x)P_1(x+u_2)Q_2(x+u_1)Q_1(x)}{P_2(x+u_1)P_1(x)Q_2(x)Q_1(x+u_2)}=\nu^n.
\]
On the left-hand side, the numerator and the denominator have the same degree and the same dominant coefficient. This implies $\nu^n=1$, which is absurd because $\nu$ has infinite order. Thus, $g_3$ is not of the form $(x,y)\mapsto (x,\nu y)$.

Hence $g_3$ is of the form $(x,y)\mapsto (x,y+v_3)$. By Theorem \ref{ellipticthm} $g_0,g_1,g_2$ can be respectively written as $(\gamma x, y+R_0(x))$ and $(x+u_i,y+R_i(x)),\ i=1,2$ where $R_0,R_1,R_2\in\CC(x)$. 

We will exploit the relation $g_1^{-1}g_2^{-1}g_1g_2=g_3^n$ to show that $R_1,R_2$ must be polynomials. Note that $g_3$ is elliptic and acts trivially on the base; roughly speaking $g_1,g_2$ almost commute. Before we continue the proof we recall first some notions. An indeterminacy point $x$ of $f$ will be called persistent if for every $i>0$, $f^{-i}$ is regular at $x$ and the backward orbit of $x$ is infinite, and if there are infinitely many curves contracted onto $x$ by the iterates $f^{-k},k\in\NN$. A conic bundle is a rational fibration where the only singular fibers are unions of two $(-1)$-curves. It is proved in \cite{Zhaocent} that $g_1$, being an element of $\Jonq$, acts by algebraically stable transformation on a conic bundle $X$; moreover the only singular fiber of $X$ lies over the $\Gamma_N$-invariant fiber $x=\infty$. 

Suppose by contradiction that $R_1$ is not a polynomial; this implies that $g_1$ is a \Jonqui twist. Some poles of $R_1$ in $\CC$ correspond to persisitent indeterminacy points of $g_1$ on $X$ (see \cite{Zhaocent} for details). Let $e\in X$ be a persistent indeterminacy point of $g_1$. Since $\{g_1^{-i}(e),i>0\}$ is infinite, $g_2$ and $g_3$ are regular at $g_1^{-k}(e)$ for $k$ large enough. For infinitely many $j>0$, $g_1^{-j}$ contracts a regular fiber of the conic bundle onto $e$, denote it by $C_j$. For $k$ large enough $g_2$ and $g_3$ do not contract $C_k$. Keeping these two observations in mind, from the relation $g_1^k\circ g_2\circ g_1^{-j}=g_2\circ g_3^{nk}\circ g_1^{k-j}$ we deduce that $g_2\circ g_1^{-j}(e)$ is an indeterminacy point of $g_1^k$ for suitable $j,k$ (recall that $g_3$ does not permute the fibers of the conic bundle). This means that, under the iteration of $g_1$, the forward orbit of $g_2\circ g_1^{-j}(e)$ will meet a persistent indeterminacy point $e'$ of $g_1$. The correspondance $e\mapsto e'$ does not depend on $j,k$. Thus, up to raplacing $g_2$ by an iterate $g_2^m$, we have $e=e'$. Then for some $l\in\ZZ$, $g_1^l\circ g_2^m( g_1^{-j}(e))$ will be an indeterminacy point of $g_1^j$, i.e.\ we have $g_1^l\circ g_2^m( g_1^{-j}(e))=g_1^{-j}(e)$. Similarly, we have $g_1^l\circ g_2^m (C_k)=C_k$ for $k$ large enough. This means that $g_1^l\circ g_2^m$ preserves the rational fibration fiber by fiber. In particular $lu_1+mu_2=0$, which is impossible because $u_1,u_2$ generate a non-discrete subgroup of $\CC$. 

Now we know that $R_1,R_2$ are polynomials. Consequently $g_1,g_2$ are elliptic. Let us finish the proof by showing that $R_0$ is a polynomial too. The element $g_0g_ig_0^{-1}$ reads 
\[
(x,y)\dashrightarrow (x+\gamma u_1,y-R_0(\gamma^{-1}x)+R_1(\gamma^{-1}x)+R_0(\gamma^{-1}x+u_1)).
\]
The relation $g_0g_ig_0^{-1}=g_1^{n_{i1}}g_2^{n_{i2}}g_3^{m_i}$ implies that the rational fraction $-R_0(\gamma^{-1}x)+R_1(\gamma^{-1}x)+R_0(\gamma^{-1}x+u_1)$ is a polynomial. This is only possible if $R_0$ is a polynomial.
\end{pro}

From the above discussions we know that $\Gamma_N$ is an elliptic subgroup of $\CRC$. The proofs of Lemma \ref{developingmapbiholomorphic} and Lemma \ref{noinvcurves} work exactly in the same way and we reduce the birational structure on $S_N$ to a $(\Aut(\FFh_k),\FFh_k)$-structure for $k\neq 1$ or to a $(\PGL _3(\CC),\PP^2)$-structure as in Section \ref{subsectioninjectiveholonomy}. If it is reduced to a $(\PGL _3(\CC),\PP^2)$-structure then the result of B. Klingler \cite{Kli98} finishes the proof. It cannot be reduced to a $(\Aut(\PP^1\times\PP^1),\PP^1\times\PP^1)$-structure because a finite unramified cover of $S_N$ would have two holomorphic foliations. 

Assume that the birational structure is reduced to a $(\Aut(\FFh_k),\FFh_k)$-structure for $k\geq 2$. Then $\Gamma_N$ preserves a rational fibration. Denote by $\sigma$ the induced homomorphism from $\Gamma_N$ to $\PGL _2(\CC)$. Using the same reasoning we have done in the proof of the previous lemma, we can write $g_0,g_1,g_2,g_3$ as:
\begin{align*}
g_0&:(x,y)\mapsto (\gamma x,y+R_0(x));\\
g_i&:(x,y)\mapsto (x+u_i,y+R_i(x)),\ i=1,2;\\
g_3&:(x,y)\mapsto (x,y+v_3)
\end{align*}
where $u_1,u_2,v_3,\gamma\in\CC^*$ and $R_1,R_2,R_3$ are polynomials. Moreover we have 
\begin{equation}
\gamma\begin{pmatrix} u_1\\ u_2\end{pmatrix}=\begin{pmatrix} n_{11}&n_{12}\\n_{21}&n_{22}\end{pmatrix}\begin{pmatrix} u_1\\ u_2\end{pmatrix}.\label{eqmatrixconj}
\end{equation}
where $\begin{pmatrix} n_{11}&n_{12}\\n_{21}&n_{22}\end{pmatrix}$ is the matrix $N$. The relation $g_1^{-1}g_2^{-1}g_1g_2=g_3^n$ reads
\begin{equation}
R_2(x)+R_1(x+u_2)-R_2(x+u_1)-R_1(x)=nv_3.
\label{eqg1g2}
\end{equation}
For the left side of Equation \eqref{eqg1g2} to be a constant, the degrees of $R_1,R_2$ must be the same. Denote by $l$ their degree. For $i=1,2$, the element $g_0g_ig_0^{-1}$ may be written as
\begin{equation}
(x,y)\mapsto (x+\gamma u_i,y-R_0(\gamma^{-1}x)+R_1(\gamma^{-1}x)+R_0(\gamma^{-1}x+u_i)).
\label{eqg0g1}
\end{equation}
The relation $g_0g_ig_0^{-1}=g_1^{n_{i1}}g_2^{n_{i2}}g_3^{m_i}$ implies that the polynomial $-R_0(\gamma^{-1}x)+R_1(\gamma^{-1}x)+R_0(\gamma^{-1}x+u_i)$ has degree $l$. This is possible only if the degree of $R_0$ is less than or equal to $(l+1)$. For $i=1,2,3$ and $0\leq j\leq l+1$, we denote by $r_{ij}$ the coefficient of $x^j$ in $R_i(x)$.

Suppose by contradiction that $l>1$. By looking at the dominant coefficients in the equations $g_1^{-1}g_2^{-1}g_1g_2=g_3^n$ and $g_0g_ig_0^{-1}=g_1^{n_{i1}}g_2^{n_{i2}}g_3^{m_i},i=1,2$, we obtain
\begin{align}
&r_{1l}lu_2-r_{2l}lu_1=0\label{eqr1r2}\\
&\gamma^{-l}r_{il}+\gamma^{-l}(l+1)u_ir_{0(l+1)}=n_{i1}r_{1l}+n_{i2}r_{2l} \quad i=1,2.\label{eqr0r1}
\end{align}
In terms of matrices, Equation \eqref{eqr0r1} reads
\begin{equation*}
(N-\gamma^{-l}\Id)\begin{pmatrix} r_{1l}\\ r_{2l}\end{pmatrix}=\gamma^{-l}(l+1)u_ir_{0(l+1)}\begin{pmatrix} 1\\ 1\end{pmatrix}
\end{equation*}
which by Equation \eqref{eqmatrixconj} and Equation \eqref{eqr1r2} is equivalent to
\begin{equation}
(\gamma-\gamma^{-l})\begin{pmatrix} u_1\\ u_2\end{pmatrix}=C\begin{pmatrix} 1\\ 1\end{pmatrix}
\label{eqmatrixprop}
\end{equation}
for some non-zero constant $C$. This is not possible because $u_1\neq u_2$. Therefore $l\leq 1$ and $g_1,g_2$ are affine transformations. The relation $g_1^{-1}g_2^{-1}g_1g_2=g_3^n$ now reads
\begin{align}
r_{11}lu_2-r_{21}lu_1=nv_3.\label{eqr1r2again}
\end{align}
Equation \eqref{eqr1r2again} implies that $l\neq 0$, i.e.\ $l=1$. Then $R_0$ is a polynomial of degree at most $2$. If $R_0$ is of degree $2$, then we can conjugate $g_0:(x,y)\mapsto (\gamma x,y+R_0(x))$ by $(x,y)\mapsto (x,y+\delta x^2)$ for an appropriate $\delta\in\CC^*$ to decrease the degree of $R_0$. Moreover the conjugation by $(x,y)\mapsto (x,y+\delta x^2)$ keeps $g_1,g_2,g_3$ affine transformations. Thus we reduce the birational structure to a complex affine structure. Using again \cite{Kli98}, we achieve the proof of Theorem \ref{mainthm} for Inoue surfaces of type $S^{\pm}$.

\nocite{}

\bibliographystyle{alpha}

\bibliography{Bibliinoue}

\vspace{8mm}
\begin{minipage}{1\textwidth}
 \begin{flushleft}
 ShengYuan Zhao\\
 Institut de Recherche Math\'ematique de Rennes\\
 Universit\'e de Rennes 1\\
 263 avenue du G\'en\'eral Leclerc, CS 74205\\
 F-35042  RENNES C\'edex\\
 \emph{e-mail:} \texttt{shengyuan.zhao@univ-rennes1.fr}\\[8mm]
 \end{flushleft}
\end{minipage}
%\hspace{0.05\textwidth}

\end{document}